\documentclass[3p]{elsarticle}

\usepackage{amsmath}
\usepackage{color}
\usepackage{hyperref}

\journal{Elsevier}

\begin{document}

\begin{frontmatter}

\title{On Meshfree GFDM Solvers for the Incompressible Navier--Stokes Equations}

\author[affil1,affil2]{Pratik Suchde\corref{mycorrespondingauthor}} 
\ead{pratik.suchde@itwm.fraunhofer.de}
\author[affil1]{J\"org Kuhnert}
\author[affil2]{Sudarshan Tiwari}

\cortext[mycorrespondingauthor]{Corresponding author}

\address[affil1]{Fraunhofer ITWM, 67663 Kaiserslautern, Germany}
\address[affil2]{Department of Mathematics, University of Kaiserslautern, 67663 Kaiserslautern, Germany}

\begin{abstract}
Meshfree solution schemes for the incompressible Navier--Stokes equations are usually based on algorithms commonly used in finite volume methods, such as projection methods, SIMPLE and PISO algorithms. However, drawbacks of these algorithms that are specific to meshfree methods have often been overlooked. In this paper, we study the drawbacks of conventionally used meshfree Generalized Finite Difference Method~(GFDM) schemes for Lagrangian incompressible Navier--Stokes equations, both operator splitting schemes and monolithic schemes. The major drawback of most of these schemes is inaccurate local approximations to the mass conservation condition. Further, we propose a new modification of a commonly used monolithic scheme that overcomes these problems and shows a better approximation for the velocity divergence condition. We then perform a numerical comparison which shows the new monolithic scheme to be more accurate than existing schemes.
\end{abstract}


\end{frontmatter}

\section{Introduction}

In hydrodynamic simulations, the computational domain is often very complex or is rapidly changing with time. For such simulation domains, the task of meshing can prove to be very challenging. Meshfree or meshless methods are becoming increasingly popular for such simulations to avoid the tasks of meshing and remeshing. These methods use the numerical basis of a set of arbitrarily distributed nodes without any underlying mesh to connect them. These nodes could either be mass carrying particles or numerical points. For each node, the only connectivity information required is a local set of neighbouring nodes over which approximations are carried out.

Meshfree Generalized Finite Difference Methods~(GFDMs) are one such class of meshfree methods. Meshfree GFDMs are strong-form methods based on weighted least squares approximations. They have been widely used (for example, \cite{Gavete2003,Katz2010,Sridar2003,Zhang2016}) and are referred under various names, including FPM, which stands for both the Finite Pointset Method \cite{Tiwari2001} and the Finite Point Method \cite{OnateFPM}; the Kinetic Meshless Method~(KMM) \cite{Praveen2007} and the Least Squares Kinetic Upwind Method~(LSKUM) \cite{Ghosh1995}. They provide several advantages over meshfree particle methods such as Smoothed Particle Hydrodynamics~(SPH) \cite{Liu2003SPH}. They provide a framework to naturally incorporate boundary conditions which is not possible in SPH. They are more adaptive in terms of the spatial discretization since they use numerical approximation points instead of mass carrying particles. These advantages, however, come at the price of a lack of strict conservation.

Several mesh-based algorithms to solve the incompressible Navier--Stokes equations have been extended to meshfree methods. These include operator splitting methods such as the projection method \cite{Chorin1968}, SIMPLE \cite{SIMPLE} and PISO \cite{PISO}. However, drawbacks of these algorithms that are specific to meshfree methods have often been overlooked. One important example of such a drawback is that the discrete Laplace operator is not the same as the discrete divergence of the discrete gradient operator. This results in inaccurate approximations to the mass conservation equation. While this problem has been solved in the context of Finite Volume Methods~(FVMs) with staggered grids, this problem persists in meshfree methods and introduces errors in most operator splitting algorithms. 

In this paper, we present a few commonly used meshfree solvers, both operator splitting solvers and monolithic solvers, for the incompressible Navier--Stokes equations and study the drawbacks of these algorithms specific to meshfree GFDMs. We then propose a new monolithic solver which attempts to overcome these drawbacks. This is done by solving an over-determined problem. The mass and momentum conservation equations are solved together, and simultaneously with a pressure-Poisson equation which is needed to improve stability. All the methods considered are spatially second order accurate and use first order time integrations. However, the new monolithic method proposed here is numerically shown to be much more accurate than existing methods. This is as a consequence of the new method providing better local approximations to the mass conservation condition.

In our earlier paper \cite{Suchde2016}, we presented a method to improve global conservative properties in meshfree GFDMs, with an application to the incompressible Navier--Stokes equations. That was done by introducing an approximate discrete divergence theorem, without changing the local accuracy of the conservation equations. In contrast, in this paper, we present a method to improve the local accuracy of the mass conservation equation.

The paper is organized as follows. In Section~\ref{sec:FPM}, we give an overview of two variants of meshfree GFDMs. In Section~\ref{sec:INSE}, we present existing meshfree GFDM algorithms used to numerically solve the incompressible Navier--Stokes equations and analyze their drawbacks. A new coupled velocity-pressure solver which addresses these drawbacks is presented in Section~\ref{sec:NEWCoupled}. Numerical results and comparisons of the new method with existing methods are presented in Section~\ref{sec:Results}, and the paper is concluded with a discussion in Section~\ref{sec:Conclusion}.

\section{Meshfree Finite Differences}
\label{sec:FPM}
In meshfree GFDMs, a cloud of $N$ numerical points is used to discretize the computational domain $\Omega$. These points, with positions $\vec{x}_i$, $i = 1, \dots, N$, are usually irregularly spaced and include both points in the interior and on the boundary $\partial\Omega$ of the domain. The only connectivity information required for each point $i$ is a set of neighbouring points $S_i$. In most cases, this neighbourhood or support $S_i$ consists of the $n=n(i)$ points closest to it, $S_i=\{\vec{x}_j : \|\vec{x}_j-\vec{x}_i\| \le h\}$, where $h=h(\vec{x},t)$ is the smoothing length, which determines the spatial discretization size. In addition to the smoothing length $h$, the point cloud is also described by the parameters $r_{max}$ and $r_{min}$. Initial set up of the point cloud, and local addition and deletion of points during a simulation is done such that there exists at least one point in every sphere of radius $r_{max}h$, and that no two points are closer than $r_{min}h$. Further details about point cloud organization and management can be found in \cite{Drumm2008,Jefferies2015,SeiboldThesis}.

There are two variations of meshfree finite differences that we consider in this paper, both of which are based on moving least squares approximations. The first is the widely used approach based on the work of Liszka and Orkisz \cite{Liszka1980}, which we refer to as the classical GFDM. This method very closely resembles traditional finite differences, and is obtained by minimizing errors obtained from Taylor expansions. The second approach is a modification of the classical GFDM based on the work of Tiwari and Kuhnert \cite{Tiwari2001} in which the errors in the PDE considered is minimized simultaneously with the errors obtained from Taylor expansions. A comparison between these approaches has been done by Illiev and Tiwari \cite{Iliev2002}. In both variations, we consider Taylor expansions up to second order terms. Higher order accuracy can be attained in both formulations and have been considered by Milewski \cite{Milewski2012}.
\subsection{Classical meshfree GFDM}
\label{sec:GFDA}

This classical method is based on the work of Liszka and Orkisz \cite{Liszka1980}. For a function $u$ defined at each numerical point $i=1,2,\dots,N$, its derivatives are approximated as
\begin{equation}
	\partial^* u(\vec x_i)\approx \tilde{\partial}^*_i u = \sum_{j\in S_i}c_{ij}^*u_j\,,
\end{equation}
where ${}^*=x,y, xx, \Delta, etc.$ represents the differential operator being approximated, $\partial^*$ represents the continuous $^*$-derivative, and $\tilde{\partial}^*_i$ represents the discrete derivative at point $i$. For a point $i$, consider Taylor expansions around it at each neighbouring point $j\in S_i$
\begin{equation}
	\label{Eq:TaylorExp}
	e_j + u(\vec{x}_j)=u(\vec{x}_i)+\nabla u \cdot (\vec{x}_j-\vec{x}_i) + \frac12 (\vec{x}_j-\vec{x}_i)^TD(\vec{x}_j-\vec{x}_i)\,.
\end{equation}
The unknown coefficients $[\nabla u, D]^T$ are computed by a weighted least squares method, by minimizing
\begin{equation}
	\label{Eq:minJ}
	\text{min } J = \sum _{j \in S_i} W_je_j^2\,,
\end{equation}
where $W$ is a weighting function used to make sure that the points closer to the central point $i$ have a larger impact than the points farther away. It is usually taken as a Gaussian distribution
\begin{equation}
	\label{Eq:WeigthingKernel}
	W_j = \exp\!\left(-\alpha \frac{\|\vec{x}_j-\vec{x}_i\|^2}{h^2}\right)\,,
\end{equation}
where $\alpha$ is a positive constant. Note that the weighting function is only defined on the local support $S_i$ consisting of $n(i)$ points. For the sake of brevity, we present only the case of one spatial dimension. Eq.\,\eqref{Eq:TaylorExp} and Eq.\,\eqref{Eq:minJ} lead to the following system which is solved at each point $i=1,\dots, N$
\begin{equation}
	\label{Eq:GFDAsystem}
	\left(\begin{array}{c}
	e_1 \\
	\vdots \\
	e_{n}
	\end{array}\right) =
		\left(\begin{array}{cc}
		\delta x_1     & \frac12\delta x_1^2    \\
		\vdots         & \vdots          \\
		\delta x_{n}   & \frac12\delta x_{n}^2  \\
		\end{array}\right)
	\left(\begin{array}{c}
	u_x \\
	u_{xx}
	\end{array}\right) -
		\left(\begin{array}{c}
		u_1 - u_i\\
		\vdots \\
		u_n - u_i
		\end{array}\right)\,.
\end{equation}
where $\delta x_j = x_j - x_i$. Or, in short form  $\vec E = M\vec a - \vec b$, with $\vec{a}=(u_x, u_{xx})^T$ and $\vec{b}=(u_1-u_i,\dots,u_n-u_i)^T$. The minimization Eq.\,\eqref{Eq:minJ} can be rewritten as
\begin{equation}
	\label{Eq:minJ2}
	\text{min } J = (M\vec a-\vec b)^TW(M\vec a - \vec b)\,,
\end{equation}
where $W$ is a diagonal matrix with entries $W_1,\dots,W_n$. A formal minimization leads to
\begin{equation}
	\label{Eq:FinalSolution_GFDM}
	\vec a = [(M^TWM)^{-1}M^TW ]\vec b\,.
\end{equation}
Which leads to the differential operator stencils
\begin{align}
	u_x\Bigr|_i &= \sum_{j \in S_i} c_{ij}^x (u_j - u_i) \,,\label{Eq:ux}\\
 	u_{xx}\Bigr|_i &= \sum_{j \in S_i} c_{ij}^{xx} (u_j - u_i) \,, \label{Eq:uxx}
\end{align}
where $c_{ij}^x$ and $c_{ij}^{xx}$ represent the values in the first and second row respectively of the matrix $[(M^TWM)^{-1}M^TW]$ in Eq.\,\eqref{Eq:FinalSolution_GFDM}. These stencils are then used to solve the PDE. For example, if we consider the PDE
\begin{equation}
	\label{Eq:examplePDE}
	au + bu_x + cu_{xx} = d\,.
\end{equation}
The derivative approximations Eq.\,\eqref{Eq:ux}, Eq.\,\eqref{Eq:uxx} are substituted into the PDE to obtain
\begin{equation}
	au_i + b\sum_{j \in S_i} c_{ij}^x (u_j - u_i) + c\sum_{j \in S_i} c_{ij}^{xx} (u_j - u_i) = d,  \qquad i=1,\dots, N
\end{equation}
which forms a large sparse implicit system which is solved with an iterative method.

An alternate, but equivalent, way to arrive at the stencil coefficients in  Eq.\,\eqref{Eq:ux} is to ensure that the derivatives of monomials $m\in\mathcal{M}$, up to the order of accuracy desired, are exactly reproduced \cite{SeiboldThesis}
\begin{align}
	\sum_{j\in S_i}c_{ij}^{x}m_j &= \partial^x_i m \qquad \forall m\in\mathcal{M}\,,\label{Eq:Consistency}\\
	\text{min } J &= \sum_{j\in S_i} W_{ij}(c_{ij}^{x})^2\,. \label{Eq:BasicMin}
\end{align}
To improve stability and to improve accuracy for Dirichlet boundary conditions, a Kronecker-delta property is often fulfilled, either by the addition of a Kronecker-delta function $\delta_{ij}$ to the monomial set~$\mathcal{M}$ \cite{Kuhnert2014} or through the weighting function $W$ \cite{Milewski2013}.
\subsection{Tiwari Approach}
\label{sec:TES}
This method is based on the work of Tiwari and Kuhnert \cite{Tiwari2001}. In this method, stencil coefficients are found that include the PDE being solved as a local constraint. Then a function is found that satisfies the conditions on all the stencils computed. The least squares procedure for minimizing the errors includes not just the Taylor expansions, but also the error in solving the PDE itself. Consider the PDE Eq.\,\eqref{Eq:examplePDE} used in the previous section. We proceed in the same way as done in Section~\ref{sec:GFDA}. The system Eq.\,\eqref{Eq:GFDAsystem} and Eq.\,\eqref{Eq:minJ} get extended to
\begin{align}
	\left(\begin{array}{c}
	e_1 \\
	\vdots \\
	e_{n} \\
	e_{n+1}
	\end{array}\right) &=
		\left(\begin{array}{ccc}
		1       & \delta x_1     & \frac12\delta x_1^2    \\
		\vdots  & \vdots         & \vdots          \\
		1       & \delta x_{n}   & \frac12\delta x_{n}^2   \\
		a       & b              &  c               \\
		\end{array}\right)
	\left(\begin{array}{c}
	u \\
	u_x \\
	u_{xx}
	\end{array}\right) -
		\left(\begin{array}{c}
		u_1 \\
		\vdots \\
		u_n   \\
		d
		\end{array}\right) \,,\label{Eq:DAsystem}\\
	\text{min } J &= \sum _{j\in S_i} W_je_j^2 + W_{n+1}e_{n+1}^2\,,	 \label{Eq:DAmin}
\end{align}
which is solved at each point $i=1,\dots, N$. The stencils to the function approximation of $u$ are the only ones of importance. Proceeding in the same way as done in Section~\ref{sec:GFDA} leads to a system similar to that obtained in Eq.\,\eqref{Eq:FinalSolution_GFDM}, the first row of which gives the function approximation stencils
\begin{equation}
	u_i = \sum_{j \in S_i} \alpha_{ij} u_j + \beta_i d \,.
\end{equation}
We then find a function $u$ that satisfies all these stencils by solving a large sparse implicit system with an iterative method.
\begin{equation}
	\label{Eq:DirectSystem}
	(1-\alpha_{ii}) u_i - \sum_{\substack{j \in S_i  \\ j \neq i  } } \alpha_{ij} u_j = \beta_i d \,.   \qquad i=1,\dots,N.
\end{equation}
This results in a diagonally dominant system. 
One of the advantages of this direct approach over the classical GFDM one is that it can be used to solve algebraically over-determined, but well-posed, systems. i.e. systems with more PDEs than variables. An important use of this is to solve PDE systems on boundary points along with the imposed boundary conditions, which is not possible by the classical GFDM. Tiwari and Kuhnert \cite{Tiwari2001} used this framework for pressure-Poisson equations in Navier--Stokes solvers. There, an over-determined system was solved on boundary points, consisting of both the boundary conditions and the pressure-Poisson equation itself. The same idea was later extended for imposing free surface boundary conditions \cite{Tiwari2002}, for compressible flows \cite{Kuhnert2003}, and recently for heat transfer boundary conditions by Res{\'e}ndiz and Saucedo \cite{Edgar2015}, among other applications. A similar idea of adding the PDE on boundaries was also done in context of the meshfree Radial Basis Functions~(RBF), by using an additional set of nodes adjacent to the boundary \cite{Fedoseyev2002}, which was shown to increase accuracy over standard RBF. Using Tiwari's approach, a similar increase in accuracy is obtained in meshfree GFDMs and can be done without the addition of extra nodes.

However, other than solving a PDE along with the boundary condition at boundary points, to the best of our knowledge, this framework has not been used to solve over-determined problems across all points in the computational domain. In this paper, we show how this framework could be used to devise a new algorithm which reduces some problems with existing meshfree GFDM solvers for the incompressible Navier--Stokes equations by solving an over-determined problem across the entire computational domain.

\section{Existing meshfree GFDM algorithms for the incompressible Navier--Stokes equations}
\label{sec:INSE}
We consider the incompressible Navier--Stokes equations in Lagrangian form.
\begin{align}
	\frac{D\vec{x}}{Dt} &=  \vec{v} \,,\label{Eq:INSE1}\\
	\nabla\cdot\vec{v}  &= 0 \,,\label{Eq:INSE2}\\
	\frac{D\vec{v}}{D t} &= \frac{\eta}{\rho}\Delta \vec{v} - \frac{1}{\rho}\nabla p + \vec{g} \,,\label{Eq:INSE3}
\end{align}
where $\vec{v}$ is the fluid velocity, $p$ is the pressure, $\rho$ is the density, $\eta$ is the dynamic viscosity and $\vec{g}$ includes both gravitational acceleration and body forces. All the algorithms considered below start with an update of point locations by solving Eq.\,\eqref{Eq:INSE1} according to
\begin{equation}
	\label{Eq:Movement}
	\vec{x}_i^{(n+1)} = \vec{x}_i^{(n)} + \vec{v}_i^{(n)}\Delta t + \frac{\vec{v}_i^{(n)} - \vec{v}_i^{(n-1)}}{\Delta t} (\Delta t)^2\,,
\end{equation}
for each point $i=1,\dots,N$, where the bracketed superscript refers to the time level. Following the movement of points, mass conservation Eq.\,\eqref{Eq:INSE2} and momentum conservation Eq.\,\eqref{Eq:INSE3} are solved according to one the methods mentioned below. 

Two broad classes of algorithms are used to solve the incompressible Navier--Stokes equations. Namely, operator splitting methods, and monolithic methods. In operator splitting methods, also referred to as fractional step methods~(FSM), partitioned methods,  or pressure-segregation methods, the momentum conservation and mass conservation equations are solved in two separate steps. In monolithic methods, also referred to as coupled velocity-pressure methods, the mass conservation and momentum conservation equations are solved together in one large system. In the remainder of this section, we present existing meshfree methods in both these classes and study their drawbacks. In the next section, we present a new monolithic meshfree method which helps reduce these problems in existing solvers.

\subsection{Operator Splitting Methods}
\label{sec:projection}

Projection methods are a type of operator splitting methods, based on the Helmholtz-Hodge decomposition \cite{Bhatia2013}, which have been widely used for approximating incompressible and weakly compressible fluid flows in Finite Volume Methods \cite{Almgren1996,Brown2001}. These method have also been commonly used in meshfree contexts (see, for example, \cite{Tiwari2002,Xu2009}). They consist of two steps. In the first step, an intermediate velocity is computed by solving the momentum conservation equation. This intermediate velocity is then projected to a divergence-free field with the help of a correction pressure to obtain the final velocity. The intermediate velocity $\vec{v}^*$ is obtained by solving
\begin{equation}
	\label{Eq:ProjectionMethods:IntermediateVelocity}
	\frac{\vec{v}^*-\vec{v}^{(n)}}{\Delta t} = \frac{\eta}{\rho}\Delta \vec{v}^* - \frac{1}{\rho}\nabla p^* + \vec{g}\,,
\end{equation}
where $p^*$ is a pressure guess. In the original projection method of Chorin \cite{Chorin1968}, the pressure guess is taken to be zero. Non-zero values of this pressure guess have been shown to produce more accurate results \cite{Brown2001}. Throughout this paper, we take the pressure guess to be the pressure at the previous time-step $p^* = p^{(n)}$. The step of computation of the intermediate velocity is often done explicitly, by replacing the $\Delta \vec{v}^*$ term in Eq.\,\eqref{Eq:ProjectionMethods:IntermediateVelocity} with $\Delta \vec{v}^{(n)} $ \cite{Tiwari2002}. However, the implicit way used in Eq.\,\eqref{Eq:ProjectionMethods:IntermediateVelocity} produces more accurate, and often more stable, results. 

In the second step, the velocity is corrected by projecting it to a divergence free space by
\begin{equation}
	\label{Eq:ProjectionMethods:ActualProjection}
	\vec{v}^{(n+1)} = \vec{v}^* - \frac{\Delta t}{\rho}\nabla p_{corr}\,.
\end{equation}
The pressure correction $p_{corr}$ is computed by a pressure-Poisson equation obtained by applying the divergence operator to Eq.\,\eqref{Eq:ProjectionMethods:ActualProjection} and setting $\nabla\cdot\vec{v}^{(n+1)}=0$
\begin{equation}
	\label{Eq:ProjectionMethods:CorretionPressure}
	\frac{\Delta t}{\rho}\Delta p_{corr} = \nabla \cdot \vec{v}^*\,.
\end{equation}
Finally the pressure is updated by
\begin{equation}
	\label{Eq:ProjectionMethods:PressureUpdate}
	p^{(n+1)} = p^* + p_{corr}\,.
\end{equation}
Different variations of such projection methods including higher order time integration have been studied, for example, by Brown \textit{et~al.}\ \cite{Brown2001}. A variety of boundary conditions have been used in these methods. A discussion on boundary conditions for projection methods in the mesh-based context can be found in Denaro\ \cite{Denaro2003}. In meshfree projection methods, a common approach is to obtain the boundary condition by projecting the underlying equation solved at interior points on the outward facing unit normal. Such approaches and the required stabilization can be found in Fang and Parriaux \cite{Fang2008} and Boroomand \textit{et~al.}\ \cite{Boroomand2005}. In this paper, we use boundary conditions dependent on the physics of the underlying problem as done in Seibold \cite{SeiboldThesis}. Further, for the spatial discretization of Eq.\,\eqref{Eq:ProjectionMethods:IntermediateVelocity} and Eq.\,\eqref{Eq:ProjectionMethods:CorretionPressure} we use the classical meshfree GFDM presented in Section~\ref{sec:GFDA}, as done, for example, by Drumm \textit{et~al.} \cite{Drumm2008}. We now consider a few drawbacks of such projection methods in the meshfree context.

\subsubsection*{Consistency of numerical differential operators: } While applying the divergence operator to Eq.\,\eqref{Eq:ProjectionMethods:ActualProjection} to obtain Eq.\,\eqref{Eq:ProjectionMethods:CorretionPressure}, the assumption is made that the Laplace operator is the same as the divergence of the gradient operator. Specifically, 
\begin{equation}
	\label{Eq:ProjectionMethods:DiffOpInconsistency}
	\nabla \cdot \nabla p_{corr} = \Delta p_{corr}\,.
\end{equation}
While this is certainly true for continuous operators, it does not hold for discrete differential operators. This leads to an error in the approximation of the numerical divergence of the new velocity. Thus, mass conservation is violated at the local level. Unlike the truncation error due to the order of spatial approximation, this error does not converge to zero with a decreasing spatial discretization.

In the context of Finite Volume Methods, this problem has been known for several decades. One proposed solution was to replace the classical discrete Laplace operator with wider stencils by taking the convolution of the divergence operator and the gradient operators, i.e. setting $\Delta := \nabla\cdot\nabla$ at the discrete level \cite{Bell1991}. Another, more widely used approach to overcome this problem is to use staggered grids in which the pressure and velocity fields are defined at different locations \cite{Harrow1965}. The staggered grid approach can not be generalized to meshfree methods since all properties, both scalars or vectors, are prescribed on the same nodes. Using the first approach of setting the numerical Laplace operator to be equal to the numerical divergence of the numerical gradient causes several issues in the meshfree context. Firstly, this would lead to different support sizes for the first and second derivatives. The second derivatives would be defined on a support double the size of that of the first derivative, which more than doubles the number of points in each support domain. Moreover, there is no control over the center stencil values, which makes diagonal dominance hard to achieve. As a result, convergence of the large linear systems can be troublesome. Thus, this problem still persists in meshfree projection methods. Not only for meshfree GFDMs, but also in various other meshfree approximation methods including SPH \cite{Cummins1999SPH,Xu2009}.

Other operator splitting algorithms such as the  PISO algorithm \cite{PISO}, SIMPLE algorithms \cite{SIMPLE} and their derivatives also rely on pressure-Poisson equations. The difference being that the Poisson equations are derived by applying the divergence operator to the momentum conservation equation. Their meshfree equivalents possess the same drawback of $\Delta \neq \nabla\cdot\nabla$ at the discrete level, which leads to the same inaccuracies in the approximation of the mass conservation condition.

\subsubsection*{Compressible boundary layer: } When using the classical GFDM approach of Section~\ref{sec:GFDA} the pressure-Poisson equation Eq.\,\eqref{Eq:ProjectionMethods:CorretionPressure} is solved at interior points with appropriate boundary conditions on the boundary points. The divergence of the velocity at boundary points depends on the pressure-Poisson equation at boundary points which is not solved at all. This results in the formation of a numerical boundary layer of compressible fluid, with non-zero divergence of velocity, during the simulation of incompressible fluids. This problem has been alleviated by Tiwari and Kuhnert \cite{Tiwari2002} by solving a system where an over-determined system is solved at boundary points, which considers both the relevant boundary conditions and the pressure-Poisson equation, by using the modified framework of Section~\ref{sec:TES}. This problem of numerical boundary layers is also avoided by several solvers that solve the pressure-Poisson system before the implicit velocity system \cite{Zhou2015}. More recently, alternate solutions to this same problem were also considered by Idelsohn and O{\~n}ate \cite{Idelsohn2010}, but specific to free surface boundaries. However, the earlier issue of $\Delta \neq \nabla\cdot\nabla$ at the discrete level is present even at the boundaries, and is a larger source of error than the compressible boundary layer.

\subsubsection*{Poor accuracy at low Reynolds flow: } Projection methods usually suffer from low accuracy, and often instability, for fluid flows at low Reynolds numbers. This problem is reduced, but not overcome, by the implicit nature of Eq.\,\eqref{Eq:ProjectionMethods:IntermediateVelocity}. As a result, projection methods do not provide good approximations for fluids with high viscosity such as molten glass \cite{Moller2007}.

\subsection{Monolithic Methods}
\label{sec:vp-}

Coupled velocity-pressure solvers usually have the advantage of being more stable, especially for larger time-step sizes, but often come at the disadvantage of ill-conditioned systems. To avoid the problem of ill-conditioned systems, Kuhnert \cite{Kuhnert2014} developed a penalty formulation based on the classical GFDM approach of Section~\ref{sec:GFDA}. In two spatial dimensions, for $\vec{v}=(u,v)$ and $\vec{g}=(g_x,g_y)$, it can be written in matrix form as
\begin{equation}
	\label{Eq:vp-}	
		\left(\begin{array}{ccc}
		I - \frac{\Delta t}{\rho}\eta\mathbf{C}^{\Delta}    && \frac{\Delta t}{\rho}\mathbf{C}^{x} \\ [5pt]
		& I - \frac{\Delta t}{\rho}\eta\mathbf{C}^{\Delta}   & \frac{\Delta t}{\rho}\mathbf{C}^{y} \\ [5pt]
		\mathbf{C}^{x} & \mathbf{C}^{y} &  -A\frac{\Delta t }{\rho}\\
		\end{array}\right)
	\left(\begin{array}{c}
	\vec{U}^{(n+1)} \\
	\vec{V}^{(n+1)} \\
	\vec{P}_{corr}
	\end{array}\right) 
	=
		\left(\begin{array}{c}
	\vec{U}^{(n)} - \frac{\Delta t}{\rho}\mathbf{C}^{x}\vec{P^*}+\Delta t \vec{G_x} \\ [5pt]
	\vec{V}^{(n)} - \frac{\Delta t}{\rho}\mathbf{C}^{y}\vec{P^*}+\Delta t \vec{G_y} \\ [5pt]
	0 \\
		\end{array}\right) \,,
\end{equation}
where $\mathbf{C}^x$ is the matrix formed by the stencil coefficients for the numerical differential operators for the $x$ derivative, $c_{ij}^x$ (see Eq.\,\eqref{Eq:FinalSolution_GFDM} and Eq.\,\eqref{Eq:ux}). Similarly for $\mathbf{C}^y$ and $\mathbf{C}^\Delta$. The vector $\vec{U}^{(n+1)}$ is formed by $u^{(n+1)}$ at all points in the computational domain, $\vec{U}^{(n+1)} = (u^{(n+1)}_1, \dots, u^{(n+1)}_N) ^T$ and similarly for the other upper case vectors. Thus, the first two blocks of rows in Eq.\,\eqref{Eq:vp-} represents the momentum conservation equation at all points. The last block of rows in Eq.\,\eqref{Eq:vp-} is a penalty formulation for the conservation of mass equation
\begin{equation}
	\label{Eq:vp-_Penalty}
	\nabla\cdot\vec{v}^{(n+1)} - A\frac{\Delta t}{\rho} \Delta p_{corr} = 0\,.
\end{equation}
Rows corresponding to boundary points in Eq.\,\eqref{Eq:vp-} are replaced with the appropriate problem specific boundary conditions. The ideal scenario would be to set $A=0$, to get an exact conservation of mass. However, that leads to an ill-conditioned system. Further, typical iterative solvers do not converge for such a system due to a complete lack of diagonal dominance in the last block of rows in Eq.\,\eqref{Eq:vp-}. While specialized solvers for saddle-point problems \cite{Benzi2005,Metsch2013}, could possibly be used to obtain solutions to these sparse linear systems, simulations are still usually unstable when $A=0$. Thus, non-zero values of $A$ need to be used. Kuhnert \cite{Kuhnert2014} takes this parameter to be in the range of $A\in(0,0.3)$. Lower values of $A$ are preferred for greater accuracy, however, higher values are needed for better conditioning and faster convergence of the resulting linear system. We note that using $A=1$ would result in an implicit and coupled projection method, with an additional velocity correction step, Eq.\,\eqref{Eq:ProjectionMethods:ActualProjection}, required. The final pressure is given as done before $p^{(n+1)} = p^* + p_{corr}$.

This penalty approach has been shown to be more stable than the meshfree projection method \cite{Jefferies2015,Kuhnert2014}, especially at low Reynolds flows. However, setting $A\neq 0$ in Eq.\,\eqref{Eq:vp-_Penalty} leads to an artificial compressibility that is numerically observed to be similar to that introduced by the inconsistency between the $\Delta$ and the $\nabla\cdot\nabla$ discrete operators in operator splitting methods. Further, this approach has the same issue of a compressible boundary layer as that in classical operator splitting methods explained earlier.

This penalty approach coupled solver and the meshfree projection method have both been widely used and have shown to be robust methods with a wide variety of applications. Under the name of the Finite Pointset Method~(FPM), they have also been used as the numerical basis of two commercially used meshfree simulation tools: NOGRID \cite{Moller2007} and the meshfree module of VPS-PAMCRASH \cite{Tramecon2013}.

\textcolor{black} Other approaches to coupled solvers include solving the momentum conservation as in the first two blocks of rows in Eq.\,\eqref{Eq:vp-}, with a pressure-Poisson equation replacing the third block of rows in Eq.\,\eqref{Eq:vp-}. Like most fractional step methods, this solves the mass conservation indirectly, and introduces the same error as mentioned earlier of $\Delta \neq \nabla\cdot\nabla$ at the discrete level. In mesh-based contexts, especially for Finite Element Methods, the equivalent systems of Eq.\,\eqref{Eq:vp-} with $A=0$ are often solved by algebraically decomposing the system which is essentially equivalent to a fractional step method \cite{Sousa2015}; or with the help of stability conditions \cite{Antunes2015} which introduce errors similar to that in the method described in this section.

\section{A New Monolithic Solver}
\label{sec:NEWCoupled}

We wish to solve the momentum and mass conservation equations as
\begin{align}
	\frac{\vec{v}^{(n+1)}-\vec{v}^{(n)}}{\Delta t} &= \frac{\eta}{\rho} \Delta \vec{v}^{(n+1)} - \frac{1}{\rho}\nabla p^* - \frac{1}{\rho}\nabla p_{corr} + \vec{g} \,,\label{Eq:ClassicalCoupled_1}\\
	\nabla\cdot\vec{v}^{(n+1)} &= 0 \,.\label{Eq:ClassicalCoupled_2}
\end{align}
We emphasize that the desire is to solve the mass conservation directly as in Eq.\,\eqref{Eq:ClassicalCoupled_2} and not indirectly via a pressure-Poisson equation. Using the classical GFDM approach of Section~\ref{sec:GFDA}, Eq.\,\eqref{Eq:ClassicalCoupled_1} and Eq.\,\eqref{Eq:ClassicalCoupled_2} lead to Eq.\,\eqref{Eq:vp-} with $A=0$. As mentioned earlier, this results in an ill-conditioned system that is very hard to solve with typical iterative procedures. This problem of ill-conditioned systems can be avoided by using Tiwari's approach presented in Section~\ref{sec:TES}. For the same, we start by rewriting Eq.\,\eqref{Eq:ClassicalCoupled_1} and Eq.\,\eqref{Eq:ClassicalCoupled_2}, in two spatial dimensions, to obtain
\begin{align}
	e_1 + u^{(n+1)} - \frac{\eta\Delta t}{\rho} \Delta u^{(n+1)} + \frac{\Delta t}{\rho}\partial^x p_{corr} &= u^{(n)} - \frac{\Delta t}{\rho}\partial^x p^* + \Delta t g_x \,,\label{Eq:NewCoupled1}\\
	e_2 + v^{(n+1)}  - \frac{\eta\Delta t}{\rho} \Delta v^{(n+1)} + \frac{\Delta t}{\rho}\partial^y p_{corr} &= v^{(n)} - \frac{\Delta t}{\rho}\partial^y p^* + \Delta t g_y \,,\label{Eq:NewCoupled2}\\
	e_3 + \nabla\cdot\vec{v}^{(n+1)} &= 0 \label{Eq:NewCoupled3}\,,
\end{align}
where $e_1$, $e_2$ and $e_3$ are the errors in the discretizations of the respective PDEs. A quadratic minimization of $e_1$, $e_2$ and $e_3$, along with the errors from the Taylor expansions in all three variables, $u$, $v$ and $p_{corr}$ can be done, similar to that done in Section~\ref{sec:TES}. This results in a system of equations similar to Eq.\,\eqref{Eq:DirectSystem} and provides better conditioned systems than Eq.\,\eqref{Eq:vp-} with $A=0$. However, resulting simulations are mostly unstable. Thus, such an approach of using Tiwari's framework of GFDM discretization to coupled solvers has not been used successfully in the past.  

A possible explanation of the instability is the lack of information about the pressure correction $p_{corr}$. To correct this and to introduce a further coupling condition between the velocity and pressure, we make use of the fact that this modified GFDM framework can be used to solve algebraically over-determined problems. The system of Eq.\,\eqref{Eq:NewCoupled1} -- Eq.\,\eqref{Eq:NewCoupled3} are augmented with a pressure-Poisson equation to form an over-determined system. The Poisson equation is obtained in a manner similar to that done by SIMPLE and PISO methods, by applying the divergence operator to the conservation of momentum equation Eq.\,\eqref{Eq:INSE3} to obtain
\begin{equation}
	\label{Eq:NewCoupled4}
	e_4 + \Delta p_{corr} + \rho \nabla\cdot \left[ \left( \vec{v}^{(n)} \cdot \nabla \right) \vec{v}^{(n+1)} \right] = \frac{\rho}{\Delta t} \nabla \cdot \vec{v}^{(n)} - \Delta p^* + \rho \nabla \cdot \vec{g} \,.
\end{equation}
Eq.\,\eqref{Eq:NewCoupled1} -- Eq.\,\eqref{Eq:NewCoupled4} are solved at all  interior points. Note that the mass conservation condition is solved directly in Eq.\,\eqref{Eq:NewCoupled3} and indirectly in Eq.\,\eqref{Eq:NewCoupled4}. At boundary points, the mass conservation condition Eq.\,\eqref{Eq:NewCoupled3} is solved in addition to the relevant boundary conditions, once again leading to more PDEs than variables. The errors in each of these PDEs or boundary conditions are minimized simultaneously with the errors in the Taylor expansions for each velocity component and the pressure. 
\begin{equation}
	\label{Eq:TES_Coupled_WeightsMin}
	\text{min } J = \sum _{j \in S_i} W_{j}(e_j^u)^2 + \sum _{j \in S_i} W_{j}(e_j^v)^2 + \sum _{j \in S_i} W_{j}(e_j^p)^2 + \sum _{k = 1}^{4} W_{\text{PDE},k}(e_{k})^2\,,
\end{equation}
where $e_j^u$, $e_j^v$ and $e_j^p$ are the errors in the Taylor expansions around point $i$ of $u$, $v$ and $p_{corr}$ respectively. For all simulations, $W_j$ are taken to as in Eq.\,\eqref{Eq:WeigthingKernel}, and $W_{\text{PDE},k}=2$. Thus, the following least squares system is solved at each interior point to obtain the function approximation stencil coefficients :  $\vec{E} = M \vec{a} - \vec{b}$ with 
\begin{align}
	\vec E &= (e_1^u,\dots\,e_n^u,e_1^v,\dots\,e_n^v,e_1^p,\dots\,e_n^p,e_1,\dots,e_4)^T\,,\label{Eq:NewCoupled_E}\\
	\vec a &= (u,u_x,u_y, u_{xx},u_{yy}, u_{xy}, v,v_x,v_y,v_{xx},v_{yy},v_{xy},p,p_x,p_y,p_{xx},p_{yy},p_{xy})^T \,,\label{Eq:NewCoupled_a}\\
	\vec b &= (u_1,\dots\,u_n,v_1,\dots\,v_n,p_1,\dots\,p_n,r_1,\dots\,r_4)^T\,,\label{Eq:NewCoupled_b}
\end{align}
where $p$ is used as shorthand for $p_{corr}$, and $r_k$ are the right hand sides of equations Eq.\,\eqref{Eq:NewCoupled1} -- Eq.\,\eqref{Eq:NewCoupled4}. Further,
%
\begin{equation}
\label{Eq:NewCoupled_M} 
M = \left(\begin{array}{ccc}
	M_T & & \\
	& M_T  & \\	
   && M_T  \\	
   M_{P1} & M_{P2} & M_{P3} \\
	\end{array}\right)\,,
\end{equation}
with $M_T$ being the parts coming from the Taylor expansions
\begin{equation}
\label{Eq:NewCoupled_M_part1} 
	M_T = \left(\begin{array}{cccccc}
1  & \delta x_1 & \delta y_1 & \frac12\delta x_1^2  & \frac12\delta y_1^2 & \delta x_1\delta y_1  \\
   & & \vdots & \vdots & &  \\
1  & \delta x_n & \delta y_n & \frac12\delta x_n^2  & \frac12\delta y_n^2 & \delta x_n\delta y_n \\
		\end{array}\right)\,,\\
\end{equation}
and the part coming from the PDEs is given by 
\begin{align}
&\left(\begin{array}{c|c|c}
M_{P1} & M_{P2} & M_{P3}
\end{array}\right) = \nonumber\\
&\left(\begin{array}{cccccc|cccccc|cccccc}
1 & 0 & 0 & -\frac{\eta\Delta t}{\rho} & -\frac{\eta\Delta t}{\rho} & 0 & 0 & 0 & 0 & 0 & 0 & 0 & 0 & \frac{\Delta t}{\rho} & 0 & 0 & 0 & 0\\
0 & 0 & 0 & 0 & 0 & 0 & 1 & 0 & 0 & -\frac{\eta\Delta t}{\rho} & -\frac{\eta\Delta t}{\rho} & 0 &  0 & 0 & \frac{\Delta t}{\rho} & 0 & 0 & 0\\
0 & 1 & 0 & 0 & 0 & 0 & 0 & 0 & 1 & 0 & 0 & 0 &  0 & 0 & 0 & 0 & 0 & 0\\
0 & \rho u^{(n)}_x & \rho v^{(n)}_x & 0 & 0 & 0 & 0 & \rho u^{(n)}_y & \rho v^{(n)}_y & 0 & 0 & 0 &  0 & 0 & 0 & 1 & 1 & 0\\
\end{array}\right) \,.\label{Eq:NewCoupled_M_part2}	
\end{align}
The first two rows in Eq.\,\eqref{Eq:NewCoupled_M_part2} represent the momentum conservation in $x$ and $y$ directions respectively, the third row represents the mass conservation equation, and the last row represents the pressure-Poisson equation. Similar systems are obtained for boundary points, with the relevant rows in Eq.\,\eqref{Eq:NewCoupled_M_part2} replaced with the boundary conditions, as done in \cite{Tiwari2002}. As done earlier, a formal minimization leads to $\vec a = [(M^TWM)^{-1}M^TW ]\vec b$. Only the function approximation stencils are of interest to us, i.e. the rows of $u$, $v$ and $p_{corr}$. They can be written as
\begin{align}
	u_i &= \sum_{j \in S_i} \alpha^u_{ij} u_j + \sum_{j \in S_i} \beta^u_{ij} v_j + \sum_{j \in S_i} \gamma^u_{ij} (p_{corr})_j + \sum_{k=1}^4\zeta^u_{ik} r_k \,,\\
	v_i &= \sum_{j \in S_i} \alpha^v_{ij} u_j + \sum_{j \in S_i} \beta^v_{ij} v_j + \sum_{j \in S_i} \gamma^v_{ij} (p_{corr})_j + \sum_{k=1}^4\zeta^v_{ik} r_k \,,\\	
	(p_{corr})_i &=  \sum_{j \in S_i} \alpha^p_{ij} u_j + \sum_{j \in S_i} \beta^p_{ij} v_j + \sum_{j \in S_i} \gamma^p_{ij} (p_{corr})_j + \sum_{k=1}^4\zeta^p_{ik} r_k \,,
\end{align}
for $i=1,\dots,N$. The coefficients $\alpha$, $\beta$, $\gamma$ and $\zeta$ represent the values in the relevant row of the matrix $[(M^TWM)^{-1}M^TW ]$. These can be rearranged to obtain the sparse linear system
\begin{align}
	(1-\alpha_{ii}^u)u_i  - \sum_{\substack{j \in S_i  \\ j \neq i  } } \alpha^u_{ij} u_j - \sum_{j \in S_i} \beta^u_{ij} v_j - \sum_{j \in S_i} \gamma^u_{ij} (p_{corr})_j &= \sum_{k=1}^4\zeta^u_{ik} r_k \,,\label{Eq:FA_1}\\	
	- \sum_{j \in S_i} \alpha^v_{ij} u_j + (1- \beta^v_{ii})v_i - 
	\sum_{\substack{j \in S_i  \\ j \neq i  }} \beta^v_{ij} v_j - \sum_{j \in S_i} \gamma^v_{ij} (p_{corr})_j &= \sum_{k=1}^4\zeta^v_{ik} r_k \,,\label{Eq:FA_2}\\	
	- \sum_{j \in S_i} \alpha^p_{ij} u_j - \sum_{j \in S_i} \beta^p_{ij} v_j + (1-\gamma_{ii}^p)(p_{corr})_i -  \sum_{\substack{j \in S_i  \\ j \neq i  }} \gamma^p_{ij} (p_{corr})_j &= \sum_{k=1}^4\zeta^p_{ik} r_k \,.\label{Eq:FA_3}
\end{align}
This resulting sparse linear system is solved using an iterative solver. Eq.\,\eqref{Eq:FA_1} -- Eq.\,\eqref{Eq:FA_3} produce a diagonally dominant system. In Eq.\,\eqref{Eq:FA_1}, for example, the magnitude of the diagonal values of $1-\alpha_{ii}^u$ are significantly larger than that of the off-diagonal values of $\alpha^u_{ij}\;j\neq i$, $\beta^u_{ij}$ and $\gamma^u_{ij}$.  
Thus, the resulting linear system converges well with typical iterative procedures.

A short comparison between the meshfree projection method presented in Section~\ref{sec:projection}, the penalty approach coupled solver of Section~\ref{sec:vp-} and the new coupled solver presented in this section is listed below

\begin{itemize}
\item All three methods have the same theoretical convergence rate with the spatial discretization as all use Taylor expansions up to the same order of accuracy. All three methods also use similar first order temporal discretizations. Higher order methods have been studied extensively, especially for mesh-based fractional step methods (for example, \cite{Zhang2014}). Such higher order approximations could be applied to all three methods considered here.

\item All three methods are ``approximate'' methods as apposed to ``exact'' methods in the sense that they only solve the mass conservation equation up to the truncation error. However, the new method does not contain the additional sources of error present in the other two methods. The projection method solves the mass conservation indirectly which leads to errors due to a lack of consistency between the first and second order derivatives. The penalty approach coupled solver attempts to solve the mass conservation directly, but introduces an artificial compressibility to improve conditioning of the resulting system. On the other hand, the new coupled solver solves the mass conservation equation directly and without introducing an artificial compressibility, and thus provides a much better approximation to the mass conservation equation than both other methods.

\item While the penalty approach coupled solver results in a compressible boundary layer, this is not present in the projection method of Tiwari and Kuhnert \cite{Tiwari2002} or the new coupled solver presented here, due to the addition of the mass balance equation on boundary points. 

\item Numerically it is observed that the new coupled solver has stability comparable to the penalty approach coupled solver, and thus, much better than the projection method.

\item Both coupled solvers solve one large implicit linear system while the projection method solves two smaller implicit linear systems. In both coupled solvers, the linear systems are of the same size, but the system is denser in the new coupled solver. The resulting system of Eq.\,\eqref{Eq:FA_1} -- Eq.\,\eqref{Eq:FA_3} can be written in matrix form as
\begin{equation}
	\label{Eq:TESCoupled_Matrix}	
		\left(\begin{array}{ccc}
		I - \boldsymbol{\alpha}^u   & -\boldsymbol{\beta}^u& -\boldsymbol{\gamma}^u \\ 
		-\boldsymbol{\alpha}^v   & I - \boldsymbol{\beta}^v& -\boldsymbol{\gamma}^v \\
		-\boldsymbol{\alpha}^p   & -\boldsymbol{\beta}^p& I - \boldsymbol{\gamma}^p \\
		\end{array}\right)
	\left(\begin{array}{c}
	\vec{U}^{(n+1)} \\
	\vec{V}^{(n+1)} \\
	\vec{P}_{corr}
	\end{array}\right) 
	=
		\left(\begin{array}{c}
	\vec{R}_1 \\
	\vec{R}_2 \\
	\vec{R}_3 \\
		\end{array}\right) \,,
\end{equation}
where $\boldsymbol{\alpha}$ , $\boldsymbol{\beta}$ and $\boldsymbol{\gamma}$ are matrices formed from the coefficients in Eq.\,\eqref{Eq:FA_1} -- Eq.\,\eqref{Eq:FA_3} and $\vec{R}_1$, $\vec{R}_2$ and $\vec{R}_3$ are vectors formed from the right hand sides of Eq.\,\eqref{Eq:FA_1} -- Eq.\,\eqref{Eq:FA_3} respectively. While the sparsity pattern of each block of rows in Eq.\,\eqref{Eq:TESCoupled_Matrix} and Eq.\,\eqref{Eq:vp-} are identical, as they are dependent only on the support domains for each point, several zero blocks are present in Eq.\,\eqref{Eq:vp-}, while all blocks are non-zero in Eq.\,\eqref{Eq:TESCoupled_Matrix}. This is no longer true for cases of spatially varying viscosity $\eta$, which are very common in fluid flow applications, for which all blocks are non-zero in the systems arising from both the coupled solvers.

\item The larger size of the implicit linear systems means that both coupled solvers have higher memory requirements than the projection method.

\item  The fastest simulation times between the three methods vary on a case by case basis, but the overall simulation time is similar for all three methods.
\end{itemize}

\section{Numerical Results}
\label{sec:Results}

The explicit time-integration for the movement of points according to Eq.\,\eqref{Eq:Movement} results in a CFL-like condition on the time-step size
\begin{equation}
	\label{Eq:VarDt}
	\Delta t= C_{\Delta t} \left( \frac{h}{\|\vec{v}\|} \right)_{min}\,.
\end{equation}
Thus, a varying time-step size according to Eq.\,\eqref{Eq:VarDt} is used in all performed simulations. To determine the numerical order of accuracy of the three methods, a small time-step is used in the first test case below by using a small value of $C_{\Delta t}$. The remaining test cases use a much larger time-step which results in lower experimental convergence rates. In the comparison of simulation times, the total time of the simulation refers to the total clock time (in seconds) of the simulation, including the initial point cloud setup and the post-processing error calculations. To ensure that these comparisons are realistic, all three methods were implemented by the same programmer and use identical memory management. All simulations were carried out in Fortran and were run serially on an Intel XeonE5-2670 CPU rated at 2.60GHz. All sparse linear systems are solved using the BiCGSTAB iterative solver \cite{BiCGSTAB} without the use of any preconditioner. For the penalty approach coupled solver, the penalty coefficient is taken to be $A=0.1$ in all simulations, since lower values result in much larger simulation times due to poor convergence of the linear systems. In all the figures below, the projection method presented in Section~\ref{sec:projection} is referred to by `Projection', the coupled solver with the penalty formulation presented in Section~\ref{sec:vp-} is referred to by `Coupled:~Penalty' and the new coupled solver which directly solves an algebraically over-determined system, as done in Section~\ref{sec:NEWCoupled}, is referred to as `Coupled:~New'. 

\subsection{Taylor-Green Vortices}
\label{sec:TaylorGreenVortices}
\begin{figure}[!htbp]
  \centering
  \includegraphics[width=0.45\textwidth]{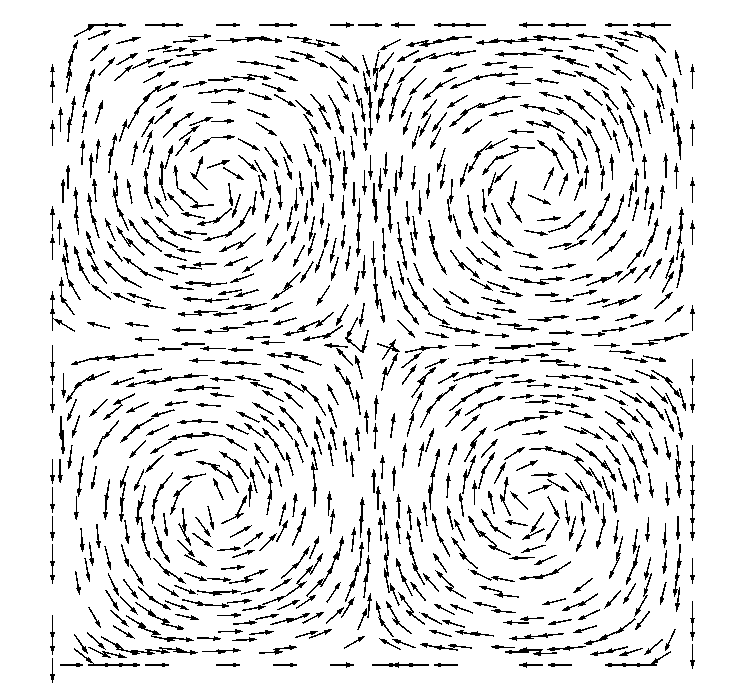}
  \caption{Initial condition for Taylor Green vortices. Arrow lengths are constant and are not scaled by velocity magnitudes.}%
  \label{Fig:TaylorGreenVortices_IC}
\end{figure}

\begin{figure}[!htbp]
  \centering
  \includegraphics[width=0.45\textwidth]{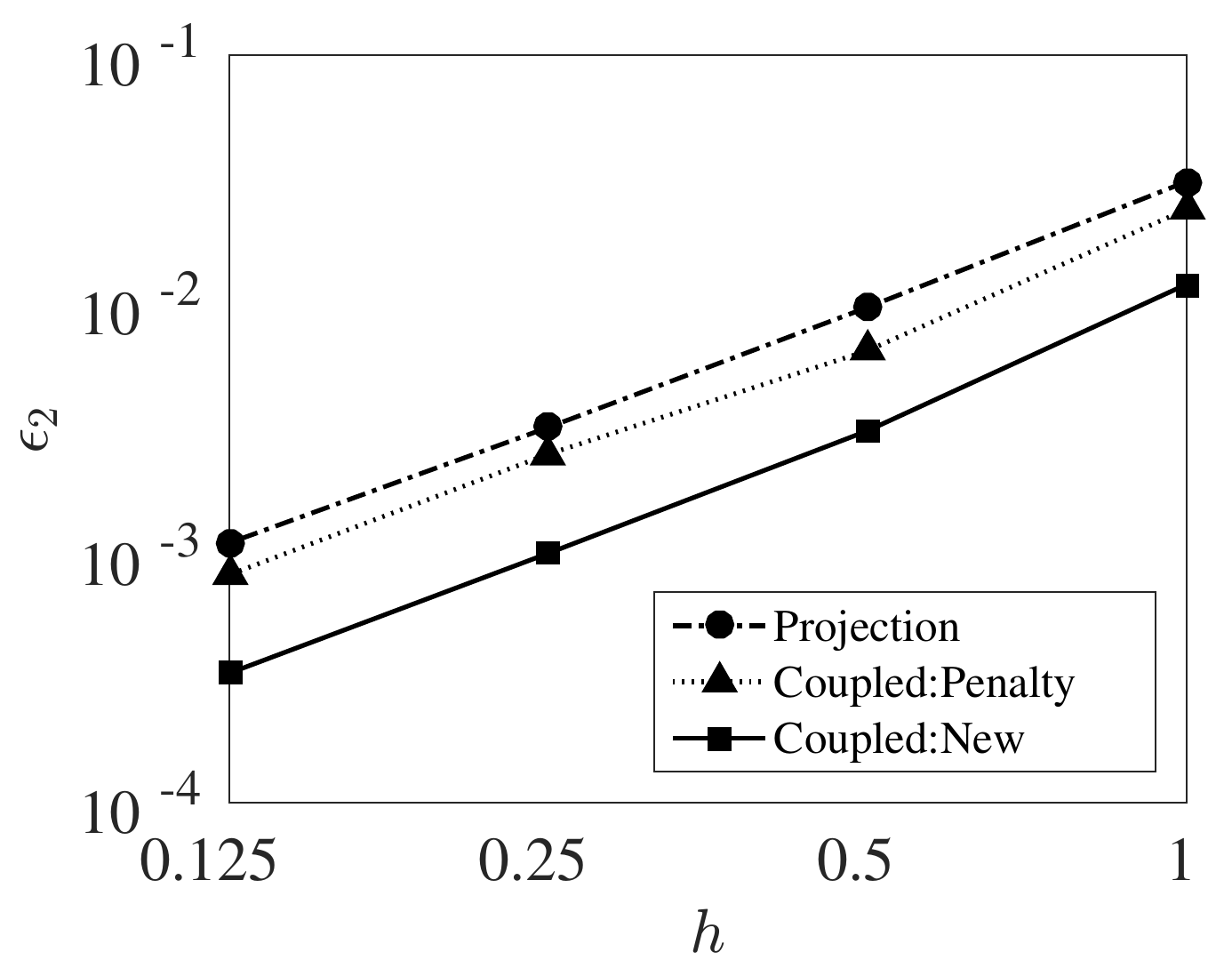}
  \caption{Convergence of error for Taylor-Green vortices.}%
  \label{Fig:TaylorGreenVortices_Results1}
\end{figure}

As a validation case, we consider the two-dimensional decaying vortices referred to as the Taylor-Green vortices on $[0,2\pi]\times[0,2\pi]$. The analytical solution is given by
\begin{align}
	u_a &= \sin(x)\cos(y)\exp(-4\pi t/Re) \,,\\
	v_a &= -\cos(x)\sin(y) \exp(-4\pi t/Re) \,,\\
	p_a & = \frac{\rho}{4}(\cos(2x)+\sin(2y))\exp(-8\pi t/Re)\,,\\
	\vec{g} & = 0\,,
\end{align}
where $\vec{v}_a = (u_a,v_a)$ is the analytical solution for the velocity, $p_a$ is the analytical solution for the pressure, and 
$Re=\frac{\rho U L}{\eta}$ is the Reynolds number, with the characteristic velocity $U=1$ and the characteristic length $L=2\pi$. Error in the numerical solution $\vec{v}$ is measured by
\begin{equation}
	 \label{Eq:PF_Error_2}
	\epsilon_2 = \left[ \frac{\sum_{i=1}^N \|\vec{v}_i - \vec{v}_{\text{a}}(\vec{x}_i)\|^2 V_i}{\sum_{i=1}^N \|\vec{v}_{\text{a}}(\vec{x}_i)\|^2V_i}	\right]^{\frac{1}{2}} \,,
\end{equation}
where $V_i$ is a volume associated with point $i$. For an error $\epsilon$ and consecutive smoothing lengths $h_1$ and $h_2$, the rate of convergence of the solution with changing smoothing length is measured as 
\begin{equation}
	r = \frac{ \log \left( \frac{\epsilon(h_2)}{\epsilon(h_1)} \right) }{ %
	\log \left( \frac{h_2}{h_1} \right)}\,.
\end{equation}
The initial condition for the velocity is taken in accordance with the exact solution, and is shown in Figure~\ref{Fig:TaylorGreenVortices_IC}. Dirichlet boundary conditions are used on all boundaries. The simulations are done up to an ending time of $t_{end} = 1s$, for $\eta = 1 Pa\,s$ and $\rho = 1 kg/m^3$, which gives $Re = 2\pi$. A small time-step is used according to Eq.\,\eqref{Eq:VarDt} with $C_{\Delta t} = 0.005$. The convergence of the errors with the smoothing length $h$ are shown in Figure~\ref{Fig:TaylorGreenVortices_Results1} and are tabulated in Table~\ref{tab:TaylorGreen}. All three methods match the analytical solution well. All three methods exhibit a similar convergence rate, which matches the theoretical expectation. Errors in the new coupled solver are smaller than the other two methods, while the errors in the projection method are the largest, but are only slightly larger than those in the penalty approach coupled solver. Total simulation times for each case are also shown in Table~\ref{tab:TaylorGreen}. All three methods take approximately the same time, with the new coupled solver being the slowest and the projection method being the fastest. Similar results were obtained for larger Reynolds numbers.

\begin{table}[!htbp]
\caption{Errors, convergence orders and simulation times for the Taylor-Green vortices test case. $h$ is the smoothing length, $N$ is the number of points in the entire domain at the initial state, $\epsilon_2$ is the relative error, $r$ is the order of convergence of $\epsilon_{2}$, and $t$ is the simulation time in seconds.}
\centering
\label{tab:TaylorGreen}
	\begin{tabular}{l|c|c|c|c|c|c|c|c|c|c}
\hline
	&& \multicolumn{3}{c|}{Projection} & \multicolumn{3}{|c|}{Coupled:Penalty}  & \multicolumn{3}{|c}{Coupled:New} \\
	$h$ & $N$ & $\epsilon_2$ & $r$ & $t$ & $\epsilon_2$ & $r$&  $t$ & $\epsilon_2$ & $r$ & $t$  \\
\hline
& &&&&&&&& \\[\dimexpr-\normalbaselineskip+2pt]
	$1.0$  & $\phantom{11\,}293$  & $3.1\times 10^{-2}$ & $-$ & $\phantom{123}6$ & $2.4\times 10^{-2}$ & $-$ & $\phantom{123}7$ & $1.2\times 10^{-2}$ & $-$ & $\phantom{123}9$\\
	$0.5$  & $\phantom{1}1\,047$ & $9.7\times 10^{-3}$ & $1.67$ & $\phantom{12}28$ & $6.5\times 10^{-3}$ & $1.88$ & $\phantom{12}29$& $3.1\times 10^{-3}$ & $1.95$ & $\phantom{12}32$\\
	$0.25$ & $\phantom{1}3\,856$ & $3.2\times 10^{-3}$ & $1.59$ & $\phantom{1}195$ & $2.5\times 10^{-3}$ & $1.37$ & $\phantom{1}202$ & $1.0\times 10^{-3}$ & $1.63$ & $\phantom{1}207$\\
	$0.125$& $14\,878$ & $1.1\times 10^{-3}$ & $1.54$ & $1934$ &  $8.2\times 10^{-4}$ & $1.60$ & $2010$ & $3.3\times 10^{-4}$ & $1.59$ & $2031$\\
\hline	
	\end{tabular}
\end{table}

\subsection{Bifurcated Tube}

\begin{figure}
	\centering
	\includegraphics[width=0.8\textwidth]{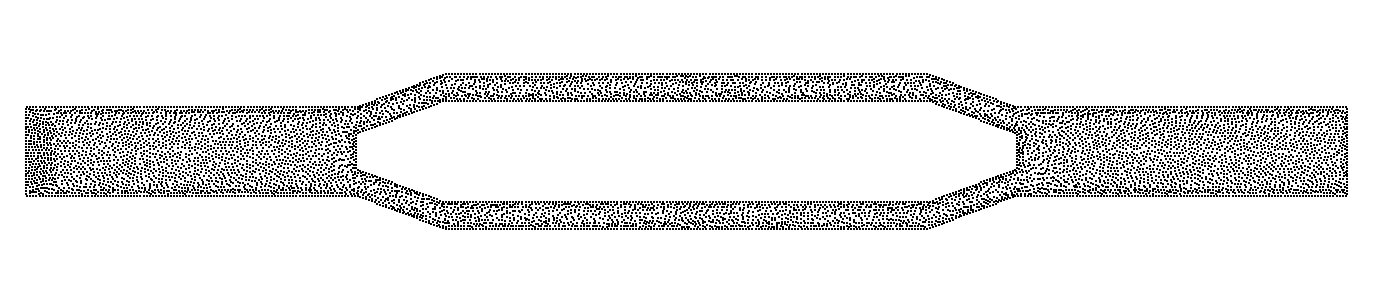}
	\caption{2D Bifurcated tube. Fluid inflow is on the left, and outflow is on the right.}
	\label{Fig:BifurcatedTube}%
\end{figure}

\begin{figure}
  \centering
  \includegraphics[width=0.45\textwidth]{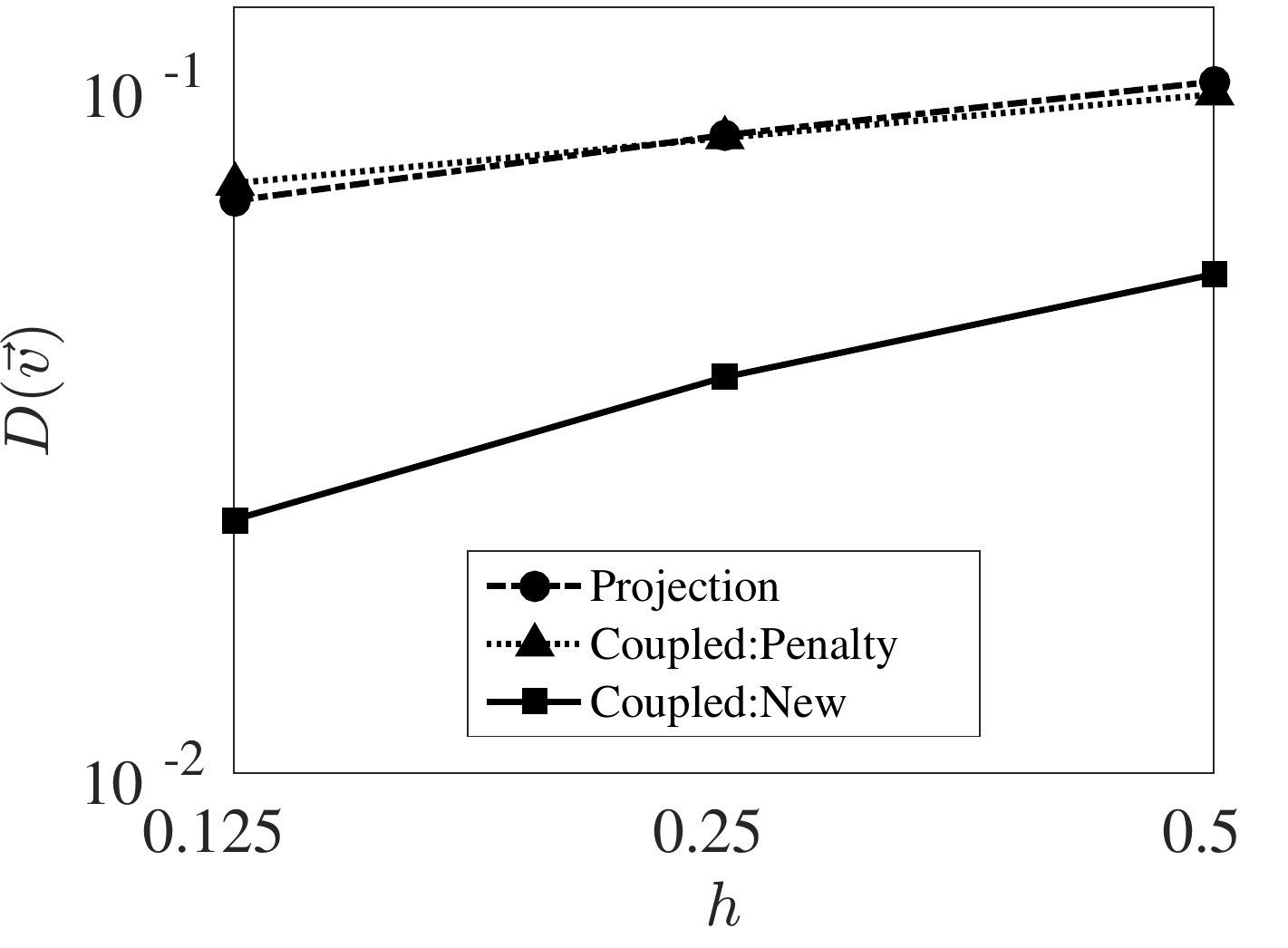}
  \includegraphics[width=0.45\textwidth]{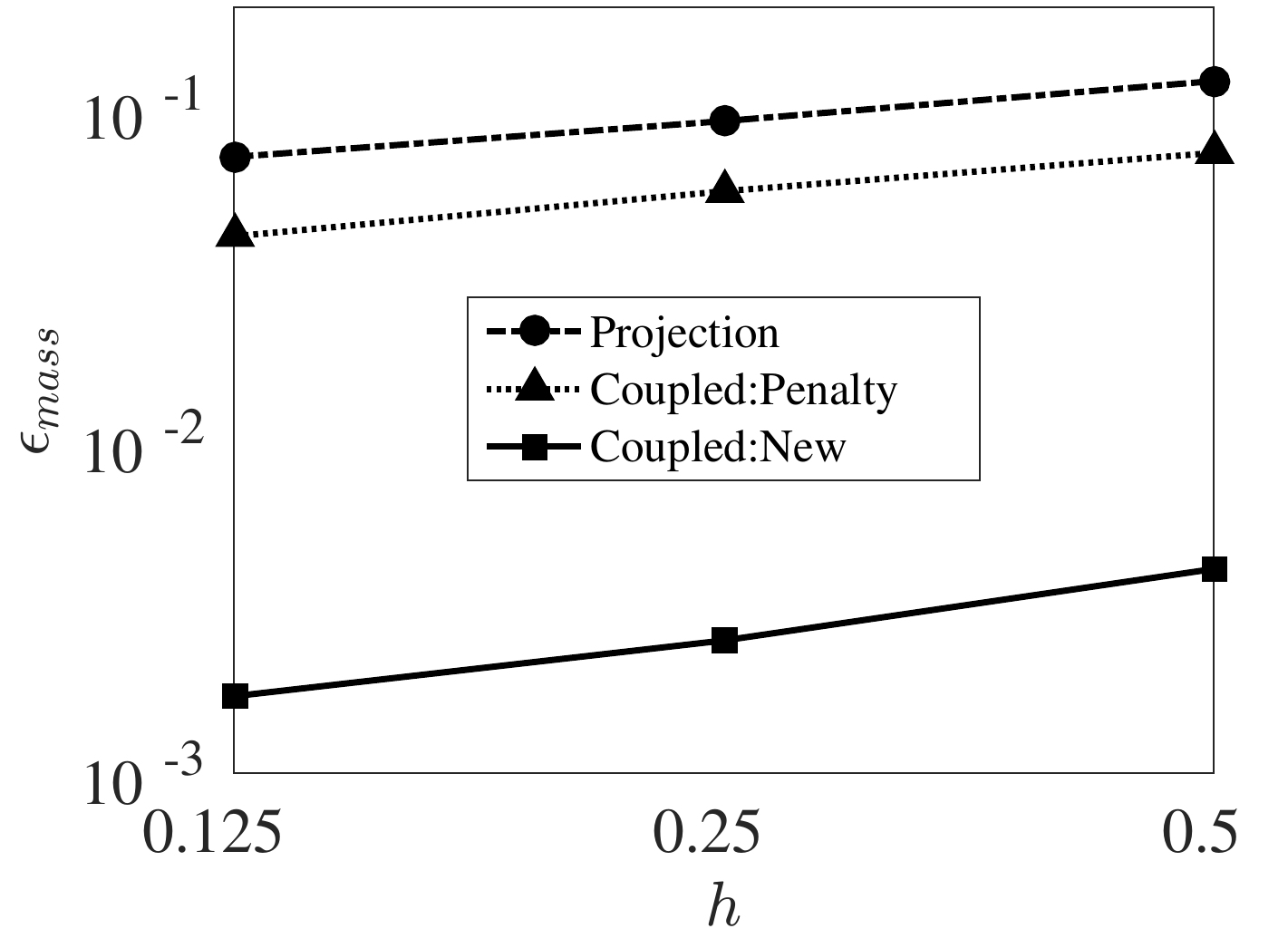}
  \caption{Bifurcated tube: average divergence of velocity in the simulation domain~(left) and error in mass conservation~(right).}%
  \label{Fig:BT_Main}
\end{figure}

We consider flow of a fluid through the bifurcated tube shown in Figure~\ref{Fig:BifurcatedTube}. The length of the tube is $60m$, the width is $4m$ in the thick region and $1m$ in the thin region.
Simulation parameters are set as $t_{end} = 1s$, $\rho = 10^3 kg/m^3$ and $\eta = 2 Pa\,s$. The velocity at the inflow, on the left of the tube, is kept constant at $\vec{v}_{in}=(2m/s,0)$. This results in a Reynolds number of the order of $10^3$. A varying time-step is used according to Eq.\,\eqref{Eq:VarDt} with $C_{\Delta t} = 0.05$. Homogeneous Neumann boundary conditions for the velocity are used at the outflow and no-slip conditions on the walls. The pressure is kept constant at atmospheric pressure at the outflow and homogeneous Neumann boundary conditions are considered elsewhere for the pressure. The error in mass conservation is measured as the difference between the total volume of fluid flowing in and that flowing out, throughout the entire simulation. The mass conservation error is given by
\begin{equation}
	\label{Eq:MassFluxError}
	\epsilon_{mass} = \left| \frac{\int_0^{t_{end}} \left[\int_{ \partial \Omega_{in}}\vec{n}\cdot\vec{v}\,dA\right] dt + \int_0^{t_{end}} \left[ \int_{\partial\Omega_{out}}\vec{n}\cdot\vec{v}\,dA\right] dt }{\int_0^{t_{end}} \left[ \int_{\partial\Omega_{in}}\vec{n}\cdot\vec{v}\,dA\right] dt} \right| \,.
\end{equation}
where $\vec{n}$ is the outward pointing unit normal and $\partial \Omega_{in}$ and $\partial \Omega_{out}$ are the inflow and outflow boundaries respectively. Note that $\epsilon_{mass}$ measures the errors during transient states too, and not just the errors in the steady state solution. A measure for the velocity divergence throughout the domain is taken as the integral of the divergence of velocity scaled by the total volume
\begin{equation}
	\label{Eq:Divergence}
	D(\vec{v}) = \frac{\int_{\Omega} |\nabla \cdot\vec{v}| dV}{\int_{\Omega} dV}\,.
\end{equation}
This can be interpreted as the average value of local error in the mass conservation equation. Note that there is no direct correlation between the measures of divergence and mass conservation, $D(\vec{v})$ and $\epsilon_{mass}$, because the absolute value of velocity divergence is taken in $D(\vec{v})$. The presence of a numerical source and a sink  of equal magnitudes for $\nabla\cdot\vec{v}$ would cancel out while measuring the mass conservation, but they would add up while measuring $D(\vec{v})$. Figure~\ref{Fig:BT_Main} shows the convergence of the velocity divergence averaged over all time-steps, and the convergence of the error in mass conservation with respect to changing smoothing length $h$. The same are also tabulated in Table~\ref{tab:BifurcatedTube}. The errors follow a similar pattern to those in the Taylor-Green vortices test case. The new coupled solver produces significantly smaller errors than the other two methods in both the average velocity divergence and mass conservation. The difference between the projection method and the coupled penalty approach is quite small. We note that convergence orders are observed to be much smaller than in the earlier test case due to the use of larger time-steps, Neumann boundary conditions and a non-standard domain geometry.
\begin{table}[!htbp]
	\caption{Errors, convergence orders and simulation times for the bifurcated tube test case. $h$ is the smoothing length, $N$ is the number of points in the entire domain at the initial state, $\epsilon_{mass}$ is the relative error in mass conservation, $r$ is the order of convergence of $\epsilon_{mass}$, and $t$ is the simulation time in seconds.}
	\centering
	\label{tab:BifurcatedTube}
	\begin{tabular}{l|c|c|c|c|c|c|c|c|c|c}
	\hline
	&& \multicolumn{3}{c|}{Projection} & \multicolumn{3}{|c|}{Coupled:Penalty}  & \multicolumn{3}{|c}{Coupled:New} \\
	$h$ & $N$ & $\epsilon_{mass}$ & $r$ & $t$ & $\epsilon_{mass}$ & $r$&  $t$ & $\epsilon_{mass}$ & $r$ & $t$  \\
	\hline
& &&&&&&&& \\[\dimexpr-\normalbaselineskip+2pt]	
	$0.5$  & $\phantom{1}5\,805$ & $1.2\times 10^{-1}$ & $-$ & $\phantom{12}556$ & $7.3\times 10^{-2}$ & $-$ & $\phantom{12}718$& $4.1\times 10^{-3}$ & $-$ & $\phantom{12}927$\\
	$0.25$ & $21\,294$ & $9.1\times 10^{-2}$ & $0.40$ & $\phantom{1}7378$ & $5.6\times 10^{-2}$ & $0.38$ & $\phantom{1}8752$ & $2.5\times 10^{-3}$ & $0.71$ & $\phantom{1}6056$\\
	$0.125$& $81\,125$ & $7.1\times 10^{-2}$ & $0.36$ & $76850$ &  $4.1\times 10^{-2}$ & $0.45$ & $76936$ & $1.7\times 10^{-3}$ & $0.56$ & $53477$\\
	\hline
	\end{tabular}
\end{table}

Simulation times for the three methods are also shown in Table~\ref{Fig:BifurcatedTube}. While the new coupled solver is the slowest by a large margin on the coarsest point cloud, it is the fastest on the finest point cloud. This could be due to slower convergence of linear systems in the other two methods.

\begin{figure}[!htbp]
  \centering
  \includegraphics[width=0.45\textwidth]{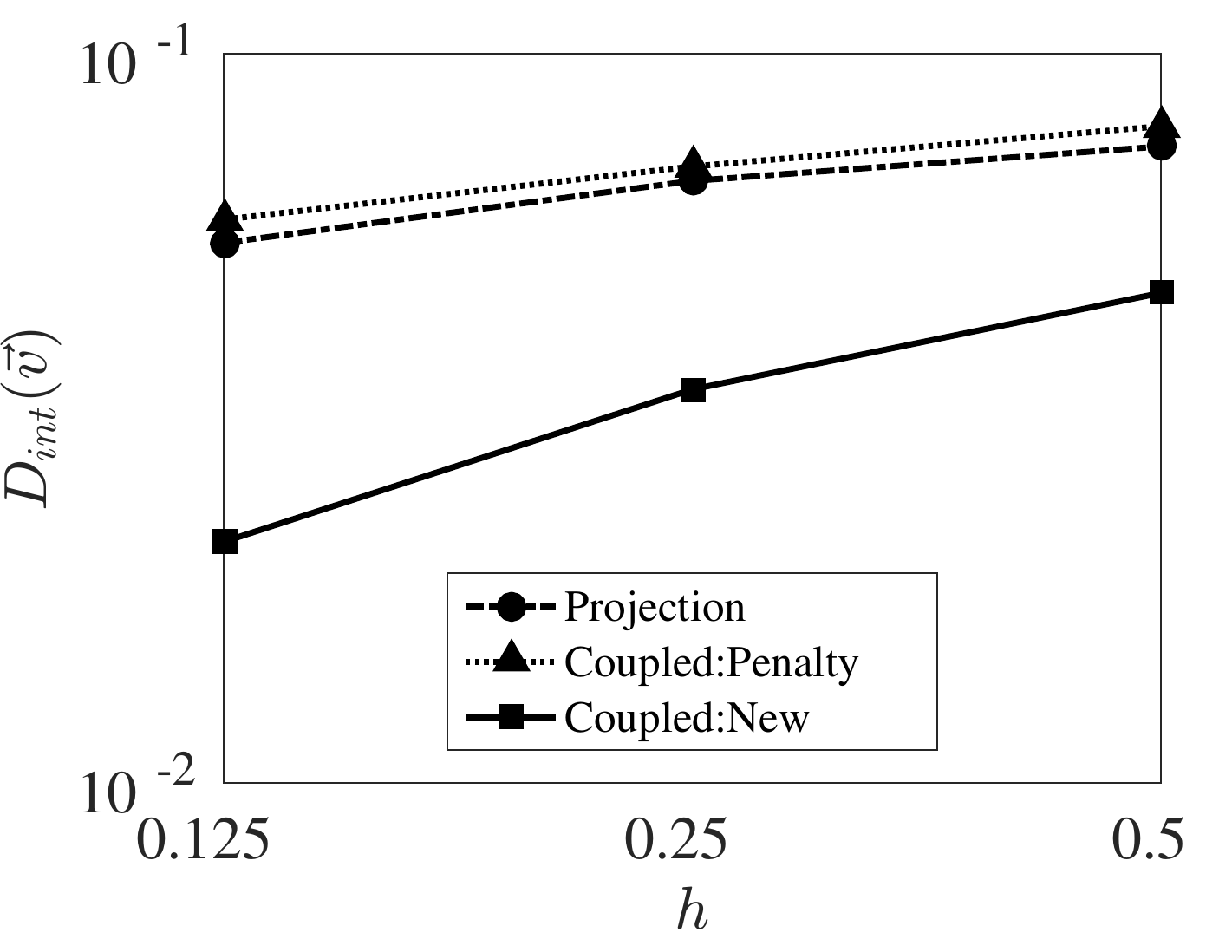}
  \includegraphics[width=0.45\textwidth]{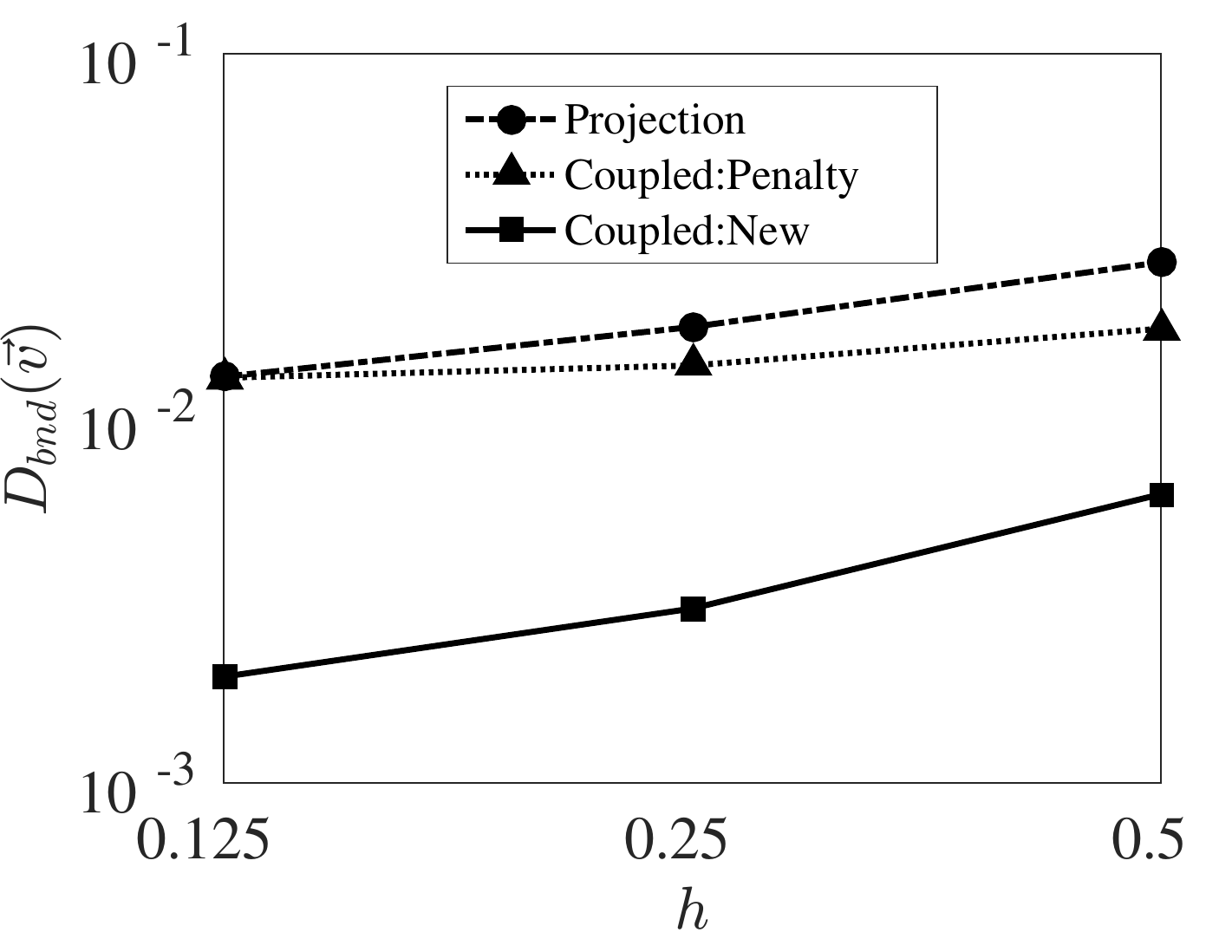}
  \caption{Velocity divergence on interior points~(left) and boundary points~(right) for the bifurcated tube test case.}%
  \label{Fig:BT_SplitDiv}
\end{figure}

To illustrate that the difference between the new coupled solver presented in this paper and the penalty formulation coupled solver goes beyond the addition of the zero divergence equation at the boundary points, we split the average divergence of velocity, Eq.\,\eqref{Eq:Divergence}, along the interior and boundary points $D(\vec{v}) = D_{int}(\vec{v}) + D_{bnd}(\vec{v})$, with
\begin{equation}
	D_{int}(\vec{v}) = \frac{\int_{\Omega\backslash\partial\Omega} |\nabla \cdot\vec{v}| dV}{\int_{\Omega} dV}\,;\;\;\; D_{bnd}(\vec{v}) = \frac{\int_{\partial\Omega} |\nabla \cdot\vec{v}| dV}{\int_{\Omega} dV}\,,
\end{equation}
where $\Omega\backslash\partial\Omega$ represents only the interior points. Figure~\ref{Fig:BT_SplitDiv} shows $D_{int}(\vec{v})$ and $D_{bnd}(\vec{v})$ averaged over all time-steps. It illustrates that the new coupled solver improves the accuracy of the mass conservation condition across both interior and boundary points.

\subsection{Sloshing}
\begin{figure}[!htbp]
  \centering
  \includegraphics[width=0.4\textwidth]{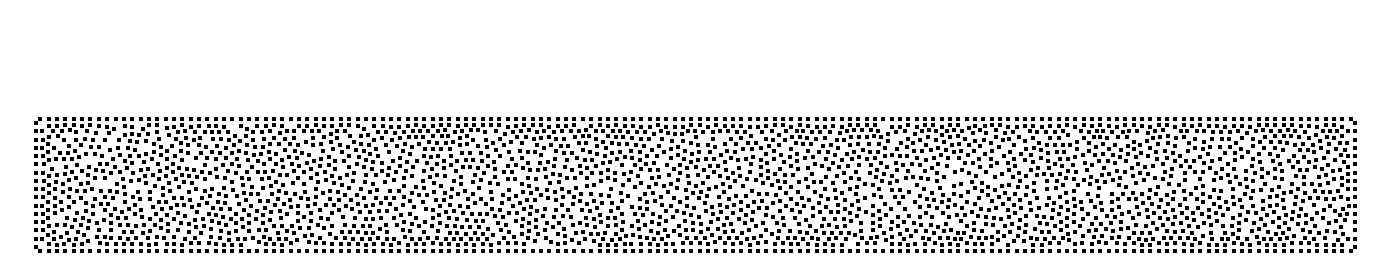}
  \includegraphics[width=0.4\textwidth]{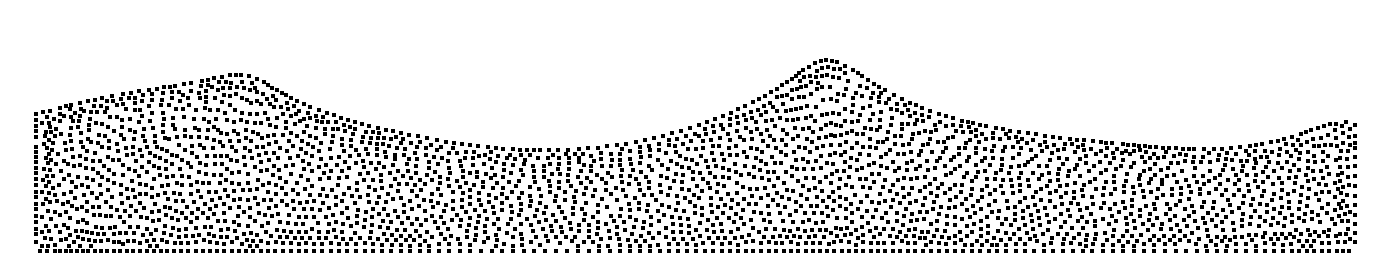}
  \vfill
  \includegraphics[width=0.4\textwidth]{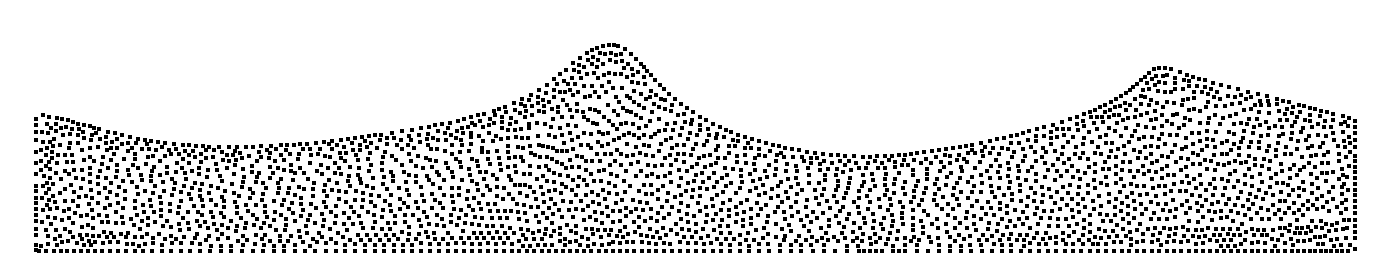}
  \includegraphics[width=0.4\textwidth]{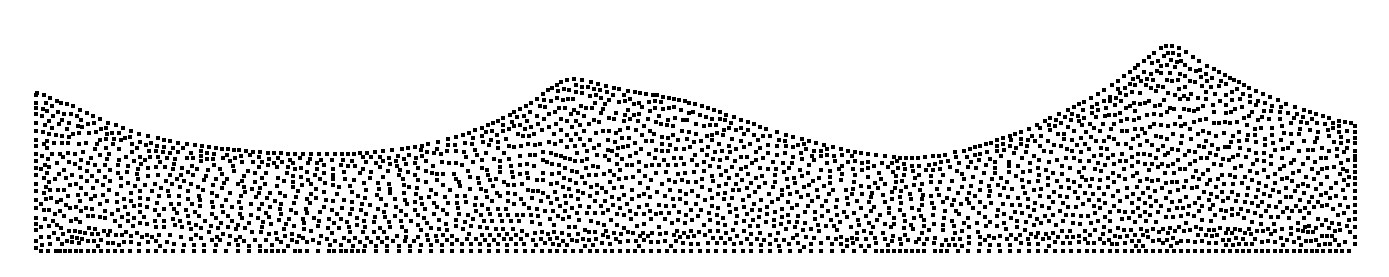}
  \caption{Sloshing at times $t = 0s$~(top left), $t = 1.31s$~(top right), $t = 1.63s$~(bottom left), and $t = 2.23s$~(bottom right).}
  \label{Fig:Sloshing}%
\end{figure}

The most common area of application of Lagrangian meshfree methods is for flows with moving free surfaces. We consider the sloshing of a fluid contained in a constantly moving rectangular box as shown in Figure~\ref{Fig:Sloshing}. The dimensions of the initial state of the fluid are $1.2m\times 0.12m$, and that of the box containing the fluid are $1.2m\times 0.6m$.

\begin{figure}
  \centering
  \includegraphics[width=0.4\textwidth]{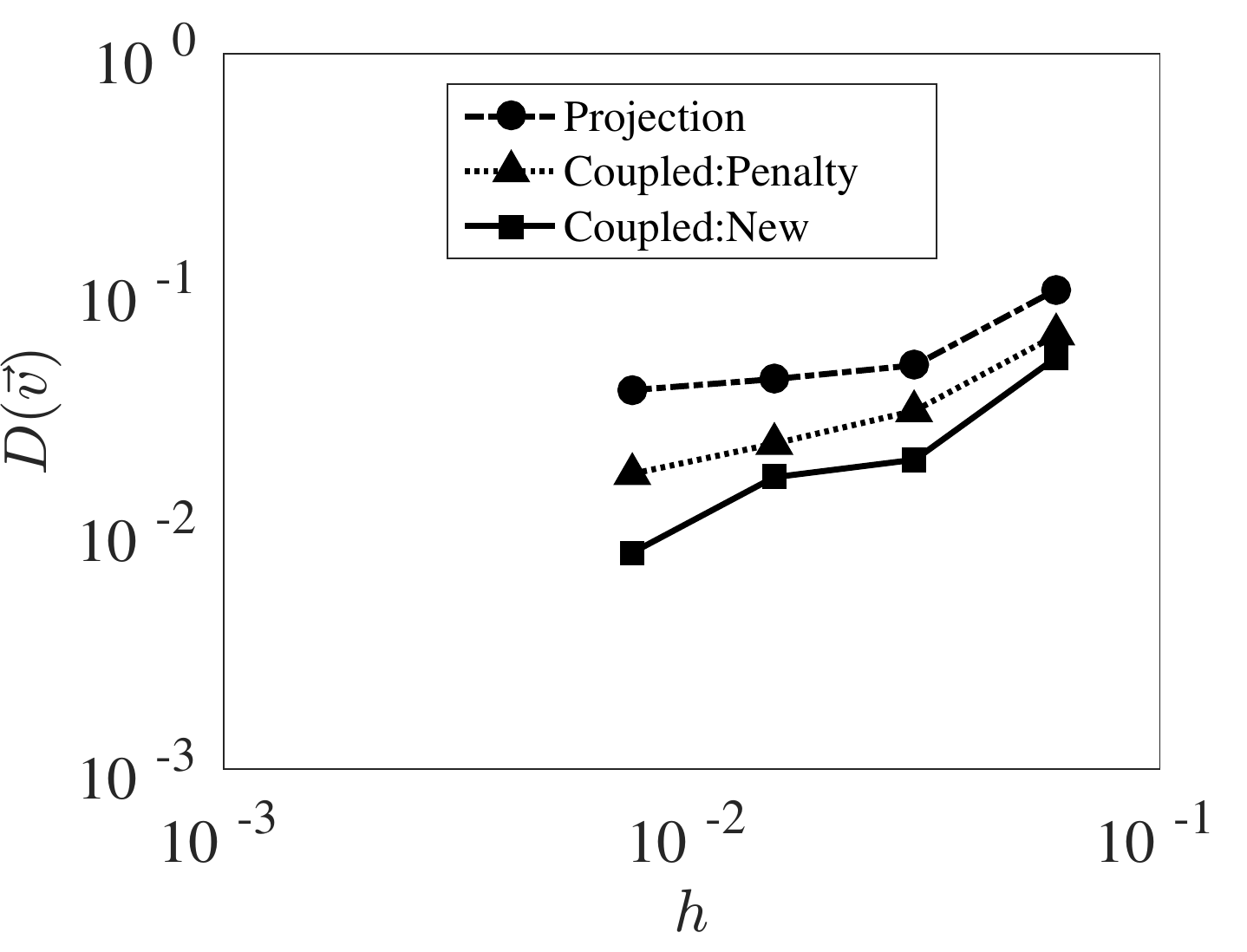}
  \includegraphics[width=0.4\textwidth]{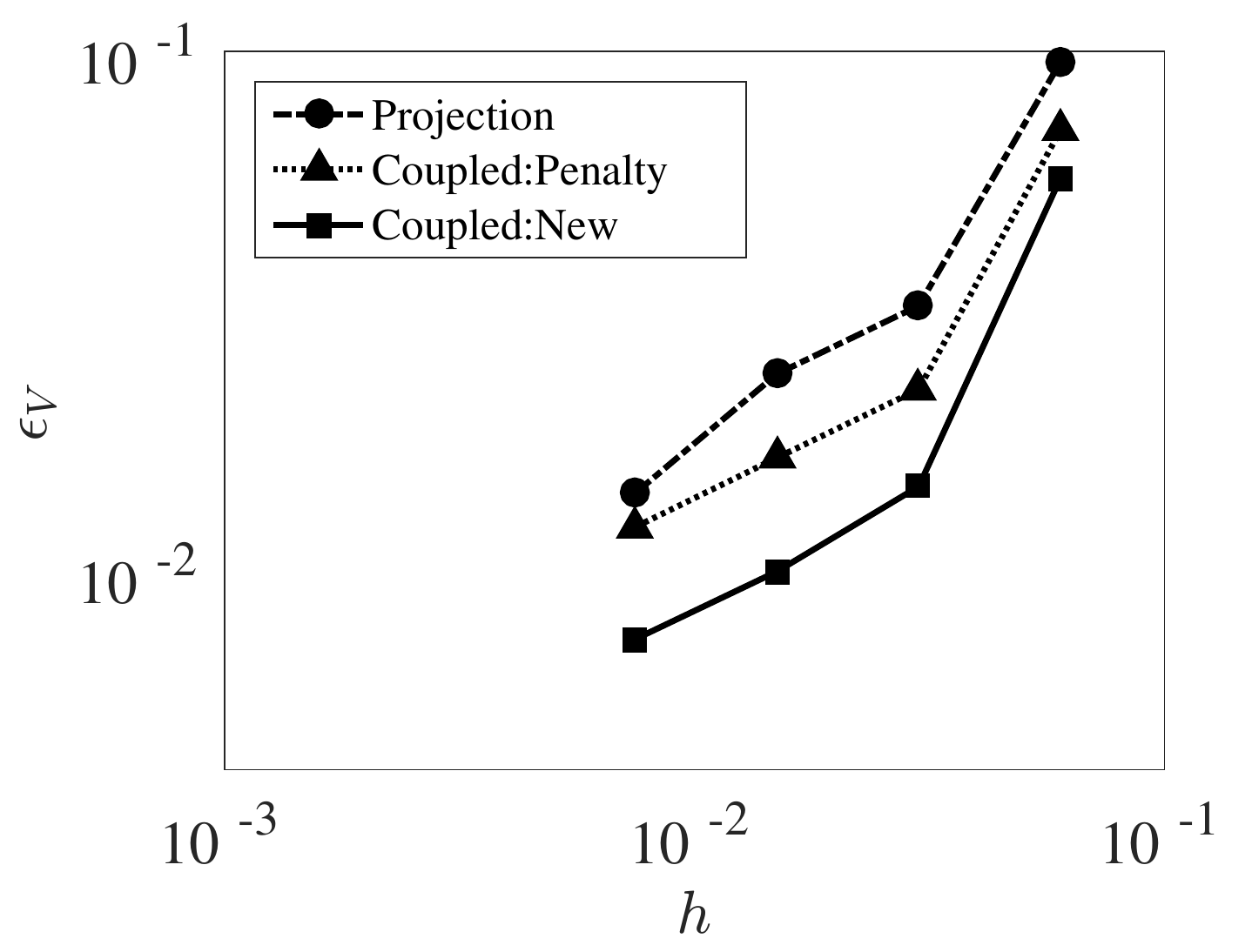}
  \caption{Sloshing: convergence of solution with smoothing length $h$. Divergence of velocity field~(left) and error in mass conservation~(right).}%
  \label{Fig:Slosh_1}
\end{figure}

The initial state is taken to be at rest. Slip boundary conditions are used at the walls for the velocity. The free surface boundary conditions are given by
\begin{align*}
	\vec{t}^T\cdot\mathbf{S}\cdot\vec{n} &= 0\\
	\vec{n}^T\cdot\mathbf{S}\cdot\vec{n} &= p - p_0 
\end{align*}
where $\vec{t}$ is a unit vector in the tangential direction to the free surface, $p_0$ is the atmospheric pressure, and $\mathbf{S}=\mathbf{S}(\vec{v})$ is the stress tensor. Homogeneous Neumann boundary conditions are used for the pressure at the walls. The movement of the box is represented in the gravitational and body forces term by setting $\vec{g}=(2\cos(10t),-10)$. The simulation parameters are set as $t_{end} = 3s$, $\rho = 10^3 kg/m^3$, and $\eta = 0.1 Pa\,s$, which leads to a Reynolds number of the order of $10^4$. A varying time-step is used according to Eq.\,\eqref{Eq:VarDt} with $C_{\Delta t} = 0.3$. The error in mass conservation is measured by the change in total volume occupied by all points, since the density $\rho$ is fixed and constant throughout the domain
\begin{equation}
	\epsilon_{V}=\frac{|\int_{\Omega_0} dV - \int_{\Omega_{end}} dV|}	{\int_{\Omega_0} dV},
\end{equation}
where $\Omega_0$ is the initial domain and $\Omega_{end}$ is the domain at $t_{end}$. The average velocity divergence is measured as done earlier in Eq.\,\eqref{Eq:Divergence}. The convergence of volume conservation error and average velocity divergence $D(\vec{v})$ with respect to a changing smoothing length $h$ for all three methods are shown in Figure~\ref{Fig:Slosh_1} and are tabulated in Table~\ref{tab:Sloshing}. The results follow a similar trend to that in the earlier two test cases. The new coupled solver shows the highest accuracy in both velocity divergence and mass conservation while the projection method shows the least accuracy in both. However, the difference between the new method and the two older methods is not as large as in the previous test case. The sharp increase in accuracy between $h=0.06$ and $h=0.03$ is due to a bad approximation of the free surface at the coarsest grid. For the finer point clouds, a small convergence rate is observed once again due to the boundary conditions used and the use of large time-steps.
\begin{table}
\caption{Errors, convergence orders and simulation times for the sloshing test case. $h$ is the smoothing length, $N$ is the number of points in the entire domain at the initial state, $\epsilon_{V}$ is the error in volume conservation, $r$ is the order of convergence of $\epsilon_{V}$, and $t$ is the simulation time in seconds.}
	\centering
	\label{tab:Sloshing}
	\begin{tabular}{l|c|c|c|c|c|c|c|c|c|c}
	\hline
	&& \multicolumn{3}{c|}{Projection} & \multicolumn{3}{|c|}{Coupled:Penalty}  & \multicolumn{3}{|c}{Coupled:New} \\
	$h$ & $N$ & $\epsilon_{V}$ & $r$ & $t$ & $\epsilon_{V}$ & $r$&  $t$ & $\epsilon_{V}$ & $r$ & $t$  \\
	\hline
	&&&&&&&&&& \\[\dimexpr-\normalbaselineskip+2pt]
	$0.06$  & $\phantom{12}\,382$ & $9.59\times 10^{-2}$ & $-$ & $\phantom{10}10$ & $7.09\times 10^{-2}$ & $-$ & $\phantom{12}12$& $5.73\times 10^{-2}$ & $-$ & $\phantom{12}15$\\
	$0.03$ & $\phantom{1}1\,361$ & $3.27\times 10^{-2}$ & $1.56$ & $\phantom{11}51$ & $2.26\times 10^{-2}$ & $1.65$ & $\phantom{12}59$ & $1.47\times 10^{-2}$ & $1.96$ & $\phantom{12}68$\\
	$0.015$& $\phantom{1}5\,119$ & $2.42\times 10^{-2}$ & $0.43$ & $\phantom{1}723$ &  $1.67\times 10^{-2}$ & $0.44$ & $\phantom{1}602$ & $1.01\times 10^{-2}$ & $0.54$ & $\phantom{1}746$\\
	$0.0075$& $19\,704$ & $1.43\times 10^{-2}$ & $0.76$ & $7871$ &  $1.23\times 10^{-2}$ & $0.44$ & $6685$ & $0.75\times 10^{-2}$ & $0.43$ & $7268$\\	
	\hline
	\end{tabular}
\end{table}

Simulation times for the three methods are also present in Table~\ref{tab:Sloshing}. All three methods take almost the same simulation time, with the new coupled solver being the slowest for 3 out of the 4 point clouds considered. 

We note that more turbulent sloshing problems than the one considered here produce larger errors in volume conservation, but for such problems, the largest source of error is in the management of the point cloud and not in the approximation of the PDEs.

\subsection{Accuracy of Gradient Reconstruction}

In the classical meshfree GFDM approach of Section~\ref{sec:GFDA}, the errors in Taylor expansions are minimized directly. As a result, first order monomials are differentiated exactly, which is evident from the formulation in Eq.\,\eqref{Eq:Consistency}. However, the same need not be true for the modified approach of Section~\ref{sec:TES} used for the new coupled solver presented in this paper. The addition of the PDE error terms to the functional minimization in Eq.\,\eqref{Eq:DAmin} and Eq.\,\eqref{Eq:TES_Coupled_WeightsMin} could result in larger errors in the Taylor expansions. This would introduce errors in gradient reconstruction.

For  $\vec{v}=(u,v)$, we examine this difference numerically by looking at the errors in the Taylor expansions for $u$. Formally, the functional being compared is $J=\sum_{j\in S_i} W_j(e_j^u)^2$ for each point $i$. After the velocity at the new time-level is computed, the obtained velocity is checked for errors in the Taylor expansions. Thus, these errors include not just the truncation error in the discretization, but also the error due to the tolerance of the sparse iterative solver. This comparison is done using the Taylor-Green vortices test case considered earlier in Section~\ref{sec:TaylorGreenVortices}, and for velocities after the completion of the first time step of the simulation. These errors are tabulated in Table~\ref{tab:GradRecon} for the classical GFDM (with the projection method) and the modified framework used for the new coupled solver. The errors are larger in the new coupled solver, but not significantly.

\begin{table}[!htbp]
	\caption{Errors in Taylor expansions measured by the functional $J=\sum_{j\in S_i} W_j(e_j^u)^2$. The table shows the mean value of $J$ across all interior points at the end of the first time-step for the Taylor-Green vortices test case.}
	\centering
	\label{tab:GradRecon}
	\begin{tabular}{|c|c|c|}
	\hline
	$h$ & Classical GFDM: Projection & Modified GFDM: New Coupled Solver  \\
	\hline
	&& \\[\dimexpr-\normalbaselineskip+2pt]
	$1.0$   & $0.2976\phantom{0}$ & $0.3057\phantom{0}$ \\
	$0.5$   & $0.09579$ & $0.09709$ \\	
	$0.25$  & $0.03255$ & $0.03345$ \\	
	$0.125$  & $0.01153$ & $0.01470$ \\	
	\hline
	\end{tabular}
\end{table}

\section{Conclusion}
\label{sec:Conclusion}

We presented a new monolithic algorithm for the incompressible Navier--Stokes equations, solved using a meshfree Generalized Finite Difference Method~(GFDM). While existing algorithms either solve the mass conservation indirectly via a pressure-Poisson equation or introduce an artificial compressibility, in the new method presented here, the mass conservation is solved directly without the introduction of any artificial compressibility. This results in improved local approximations for the mass conservation equation for both interior and boundary points, which results in better accuracy overall. In this method, the momentum and mass conservation equations are solved together, and simultaneously with a pressure-Poisson equation. The addition of this Poisson equation is essential for stabilizing the scheme and results in an over-determined system of PDEs. Accuracy is further improved at boundary points by solving the mass conservation equation in addition to the usual boundary conditions.

This new scheme was compared with two existing methods: one fractional step method and one monolithic method. All three methods have the same spatial and temporal order of accuracy. Numerical comparisons were done for different Reynolds flows, and the new method was shown to produce more accurate results. 

Our future work in this direction will be devoted to extending this algorithm to three spatial dimensions. A further interesting point of study is the impact of the ratio of weights in the minimization of the functional in Eq.\,\eqref{Eq:TES_Coupled_WeightsMin}. Higher weights could be used, for example, for the minimization of the error in the mass conservation equation. 

\end{document}